\documentclass{article}
\usepackage[utf8]{inputenc}
\usepackage[margin=1in, includehead, includefoot]{geometry}
\usepackage{amsthm,amssymb,amsmath}
\usepackage{enumitem}
\usepackage{mathrsfs}
\usepackage{mathtools}
\usepackage{fancyhdr}
\usepackage{indentfirst}
\usepackage{graphicx}
\usepackage{subfigure}
\usepackage{placeins}
\usepackage{wrapfig}
\usepackage{float}

\usepackage{xcolor}

\usepackage{lineno}
\newcommand{\toroman}[1]{\textit{\expandafter{\romannumeral #1\relax}}}

\title{Graphical parameters for classes of tumbling block graphs}
\author{{Suk J. Seo}  \\
{\small Computer Science Department } \\
{\small  Middle Tennessee State University } \\
{\small Murfreesboro, TN 37132,  U.S.A } \\
{\small sseo@mtsu.edu } \\ \\
\and
{Peter J. Slater}  \\
{\small Mathematical Sciences Department } \\
{\small and Computer Science Department } \\
{\small  University of Alabama in Huntsville } \\
{\small Huntsville, AL 35899 U.S.A } \\
}
\date{}

\begin{document}
\maketitle
\thispagestyle{empty}

\begin{abstract}
 The infinite tumbling block graph is a bipartite graph, where each vertex in one partite set is of degree 3 and each vertex in the other partite set is of degree 6. It is a 2-dimensional array of blocks of seven vertices and nine edges, a planar graph that has 3-D looks. This paper introduces tumbling block graphs and considers various graphical parameters for different classes of infinite and finite tumbling blocks.
\end{abstract}



\bigskip
\section{Introduction}
The $k$-cube $Q_k$ on $n = 2^k$ vertices is important for its applications in coding theory because its vertices represent the set of binary $k$-tuples.
 Further, because many graph parametric values are known for $Q_k$, graph
 heuristics can be tested and compared on these graphs.  See, for example, \cite{13} comparing heuristics to find the maximum independence number.  Also,
 much work has been done to determine values for parameters such as
 the domination number, the locating-dominating number, identifying-code number, or
 open-locating-dominating number for other graphs such as the infinite planar
 square grid, hexagonal grid, and triangular grid.
 (See, for example, [1-6, 8, 9, 11, 13, 14, 15, 17, 19-21].)
 In this paper we introduce the study of the finite and infinite planar tumbling-block graphs, classes of graphs which provide interesting examples for testing heuristics.

\begin{figure}[ht]
\centering
\includegraphics [width = 210pt]{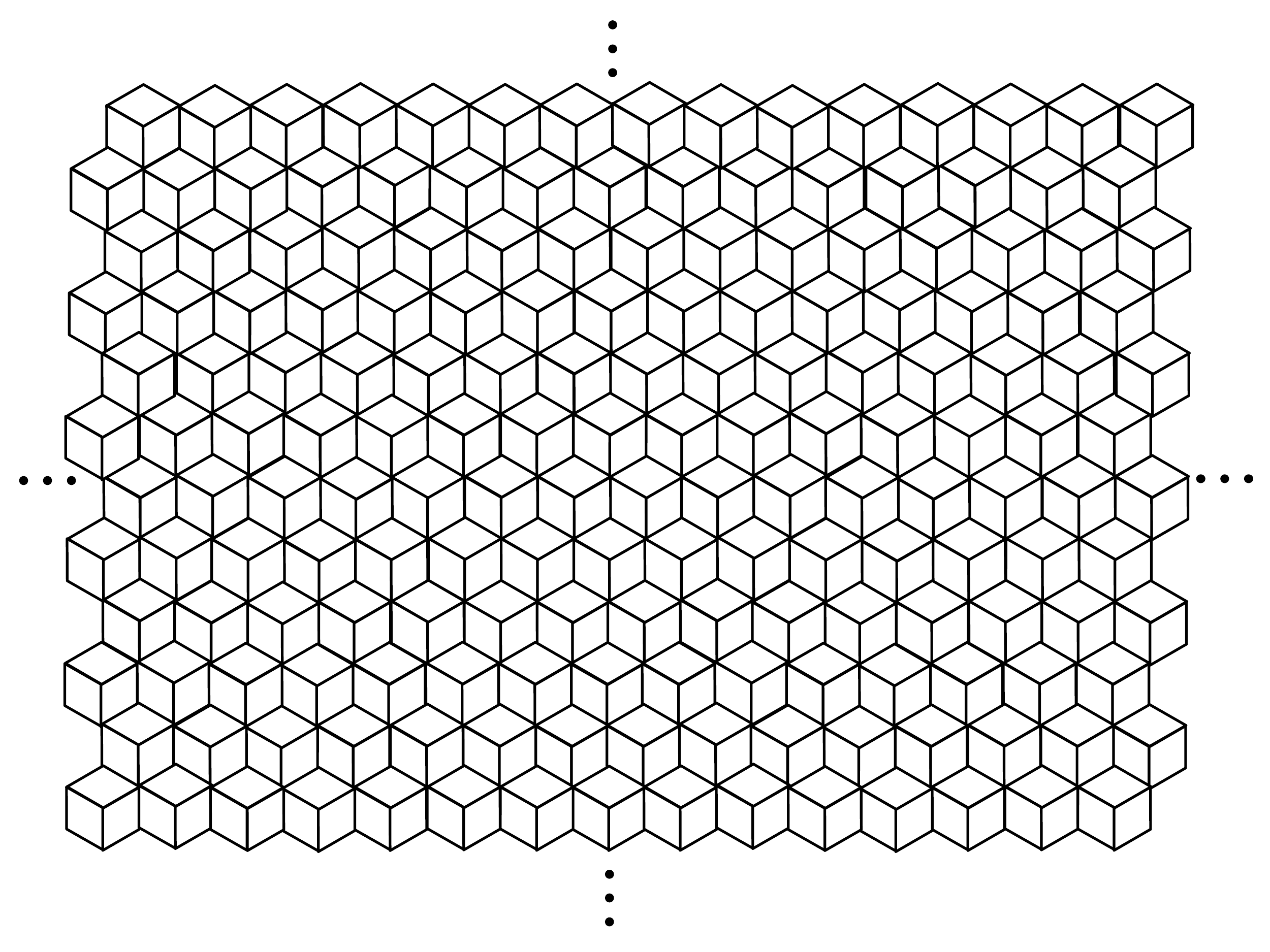}
\caption{The infinite tumbling block ($TB$)}\label{fig:tb}
\end{figure}

 The infinite tumbling block graph denoted by $TB = (V, E)$ is a two dimensional array of block graphs each of which
contains seven vertices and nine edges as shown in Figure~\ref{fig:tb}.
It is a planar graph, but it has 3-D looks. We observe that each block is $Q_3 - \{v\}$,
where $Q_3$ is the 3-cube and $v \in Q_3$.  It can also be viewed as a $C_6$ (or hexagon) subdivided into
three ${C_4}'s$ (or diamonds).  Each block shares each of its six outer vertices with a different set of two other blocks.
The infinite tumbling block graph is a bipartite graph with $V = V_1 \cup V_2$ and $V_1 \cap V_2 = \emptyset$,
where each vertex in $V_1$ is of degree 3 and each vertex in $V_2$ is of degree 6.
The density (that is, overall fraction of vertices contained therein) of each of $V_1$ and $V_2$ are 2/3 and 1/3, respectively.

\FloatBarrier

\begin{figure}[ht]
\centering
\includegraphics [width = 406pt]{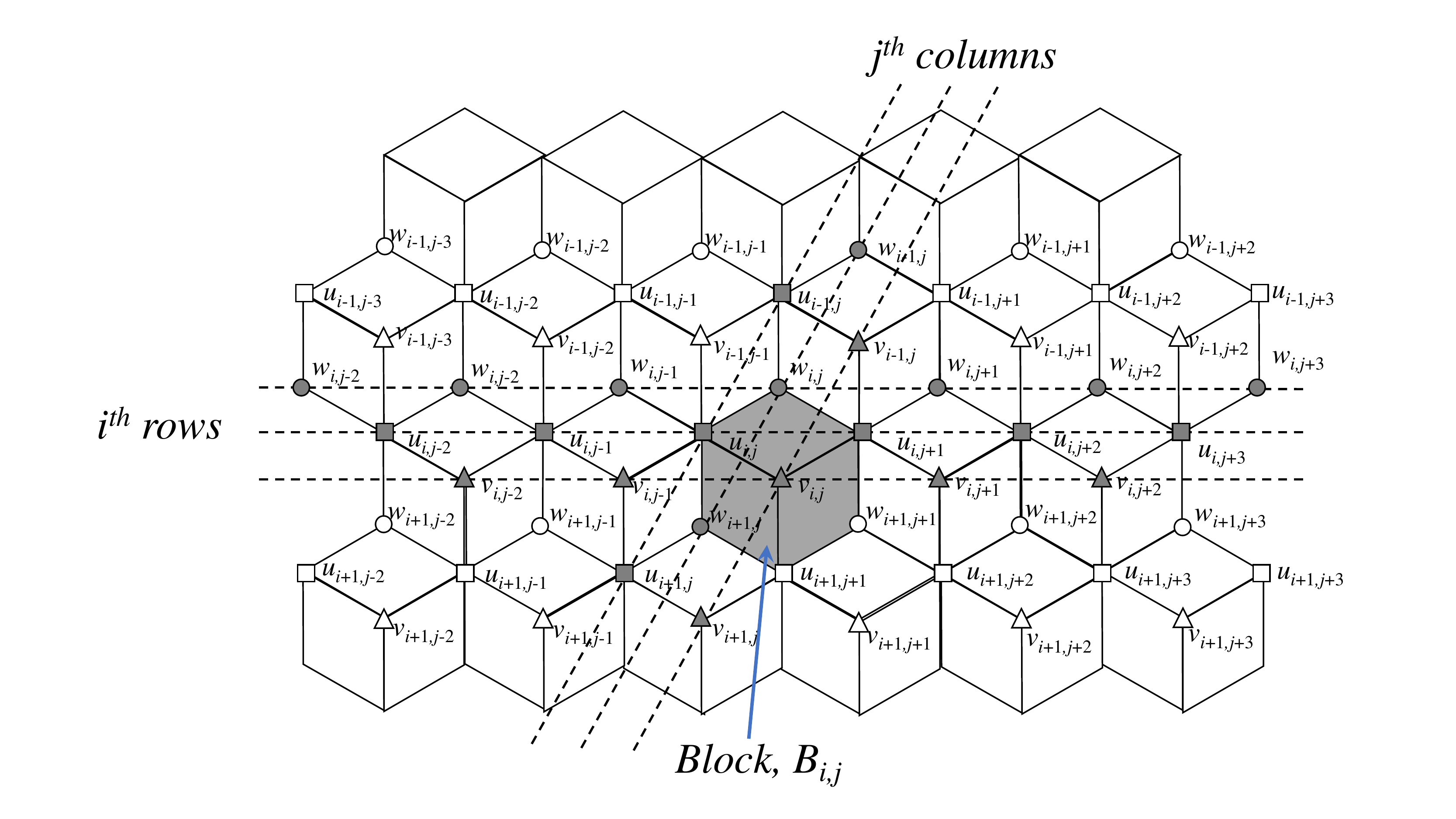}
\caption{Labeling of the vertices of an infinite tumbling block}\label{fig:tb-block}
\end{figure}

Each block of $TB$ is denoted by $B_{i,j}$, where $i$ and $j$ are the row and the column indices, respectively. We label the vertices of $TB$ as shown in Figure~\ref{fig:tb-block}.  Each block is identified by its top vertex (denoted by $w$ vertex), top-left (denoted by $u$ vertex), and middle vertex (denoted by $v$ vertex). For example, $B_{i,j}$ is identified by $w_{i,j}$, $u_{i,j}$, and $v_{i,j}$. Vertex sets $W_i^r = \{ . . . w_{i,j-2}, w_{i,j-1},w_{i,j}, w_{i,j+1}, w_{i,j+2}, . . . \}$, $U_i^r = \{ . . . u_{i,j-2}, u_{i,j-1},u_{i,j}, u_{i,j+1}, u_{i,j+2}, . . . \}$, and $V_i^r$ = ${\{ . . . v_{i,j-2}, v_{i,j-1},v_{i,j}, v_{i,j+1}, v_{i,j+2}, . . . \}}$ form the $i^{th}$ rows of $TB$.  Likewise, $W_j^c$ = $\{ . . . w_{i-2,j}, w_{i-1,j},w_{i,j}, w_{i+1,j}, w_{i+2,j}, . . . \}$, $U_j^c = \{ . . . u_{i-2,j}, u_{i-1,j},u_{i,j}, u_{i+1,j}, u_{i+2,j}, . . . \}$, and $V_j^c = \{ . . . v_{i-2,j}, v_{i-1,j},v_{i,j}, v_{i+1,j}, v_{i+2,j}, . . . \}$ form the set of $j^{th}$ columns of $TB$ as shown in Figure \ref{fig:tb-block}.

\FloatBarrier

\bigskip
Several interesting families of finite tumbling block graphs can be defined, as well as
some infinite proper subgraphs of $TB$.
We consider three families of tumbling block graphs (which have both finite and
infinite forms) based on their shapes:
\vspace{5mm}

$\bullet$ Tumbling block triangle ($TBT$)

$\bullet$ Tumbling block parallelogram ($TBP$)

$\bullet$ Tumbling block rectangle ($TBR$)

\vspace{4mm}

\begin{figure}
\centering
\subfigure[] % caption for subfigure a
{
    \label{fig:sub:a-tbt}
    \includegraphics[width=180pt]{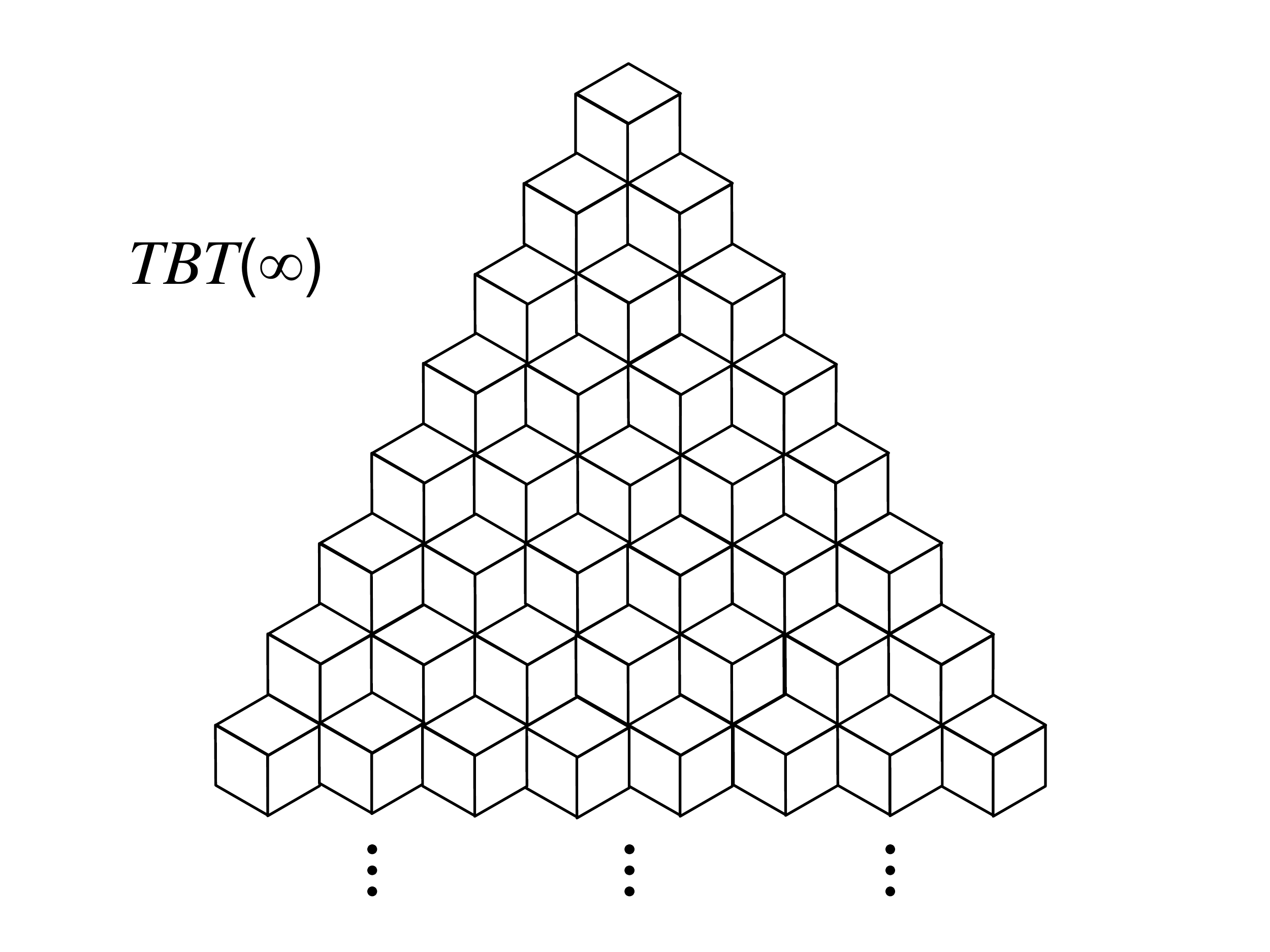}
}
\subfigure[] % caption for subfigure b
{
    \label{fig:sub:b-tbt}
    \includegraphics[width=180pt]{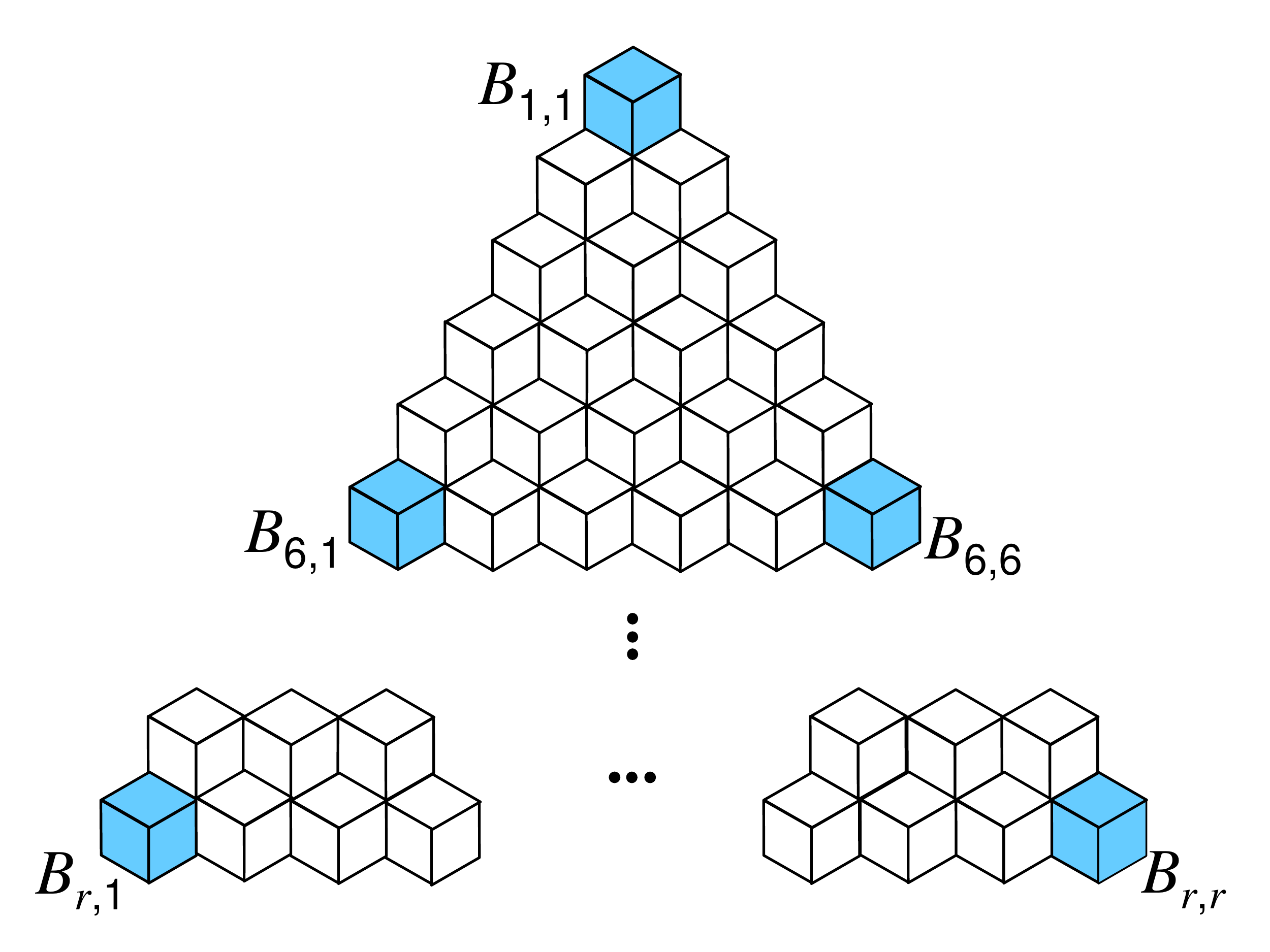}
} \caption{$TBT(\infty)$, $TBT(6)$, and $TBT(r)$ }
\label{fig:tbt} % caption for the whole figure
\end{figure}

Figure \ref{fig:tbt} shows some examples of tumbling block triangles: (a) the infinite tumbling block triangle, $TBT(\infty)$, and (b) the tumbling block triangle with height 6, $TBT(6)$,  within the tumbling block triangle with height $r$, $TBT(r)$.

\FloatBarrier

We observe that the number of vertices in the top block of $TBT(r)$,  $|V(B_{1,1})|$ is 7, the blocks in the second row add (5 + 4) additional vertices, in row three
(5 + 3 + 4) additional vertices, and so on. Hence, the number of the vertices in the tumbling block triangle with height $r$ is $|V(TBT(r))| = 7 + (5+4) + (5 + 3 +4) + ... + (5 + 3+ 3+ ... +3 + 4) = 7 + (5+4)(r-1) + 3(1+2+ ... + (r-2))
= 7 + 9r - 9 + 3(r-2)(r-1)/2 = (3r^2 +9r +2)/2$.  
Similarly, the number of edges in the tumbling block triangle with height $r$, $|E(TBT(r))| = 9 + (8+7) + (8 + 6 +7) + ... (8 + 6+ 6+ ... +6 + 7) = 9 + (8+7)(r-1) + 6(1+2+ ... + (r-2)) 	= 9 + 15r - 15 + 6(r-2)(r-1)/2	= 3r^2 + 6r$.
		
% \begin{figure}
% \centering
% \subfigure[] % caption for subfigure a
% {
%     \label{fig:sub:a-par}
%     \includegraphics[width=140pt]{fig/Figure4a.pdf}
% }
% \subfigure[] % caption for subfigure b
% {
%     \label{fig:sub:b-par}
%     \includegraphics[width=140pt]{fig/Figure4b.pdf}
% }
% \subfigure[] % caption for subfigure c
% {
%     \label{fig:sub:c-par}
%     \includegraphics[width=140pt]{fig/Figure4c.pdf}
% }
% \subfigure[] % caption for subfigure d
% {
%     \label{fig:sub:d-par}
%     \includegraphics[width=140pt]{fig/Figure4d.pdf}
% } \caption{Various tumbling block parallelograms}
% \label{fig:tbp} % caption for the whole figure
% \end{figure}
% \FloatBarrier

\begin{figure}[h]
    \centering
    \begin{tabular}{c@{\hskip 5em}c}
         \includegraphics[width=0.35\textwidth]{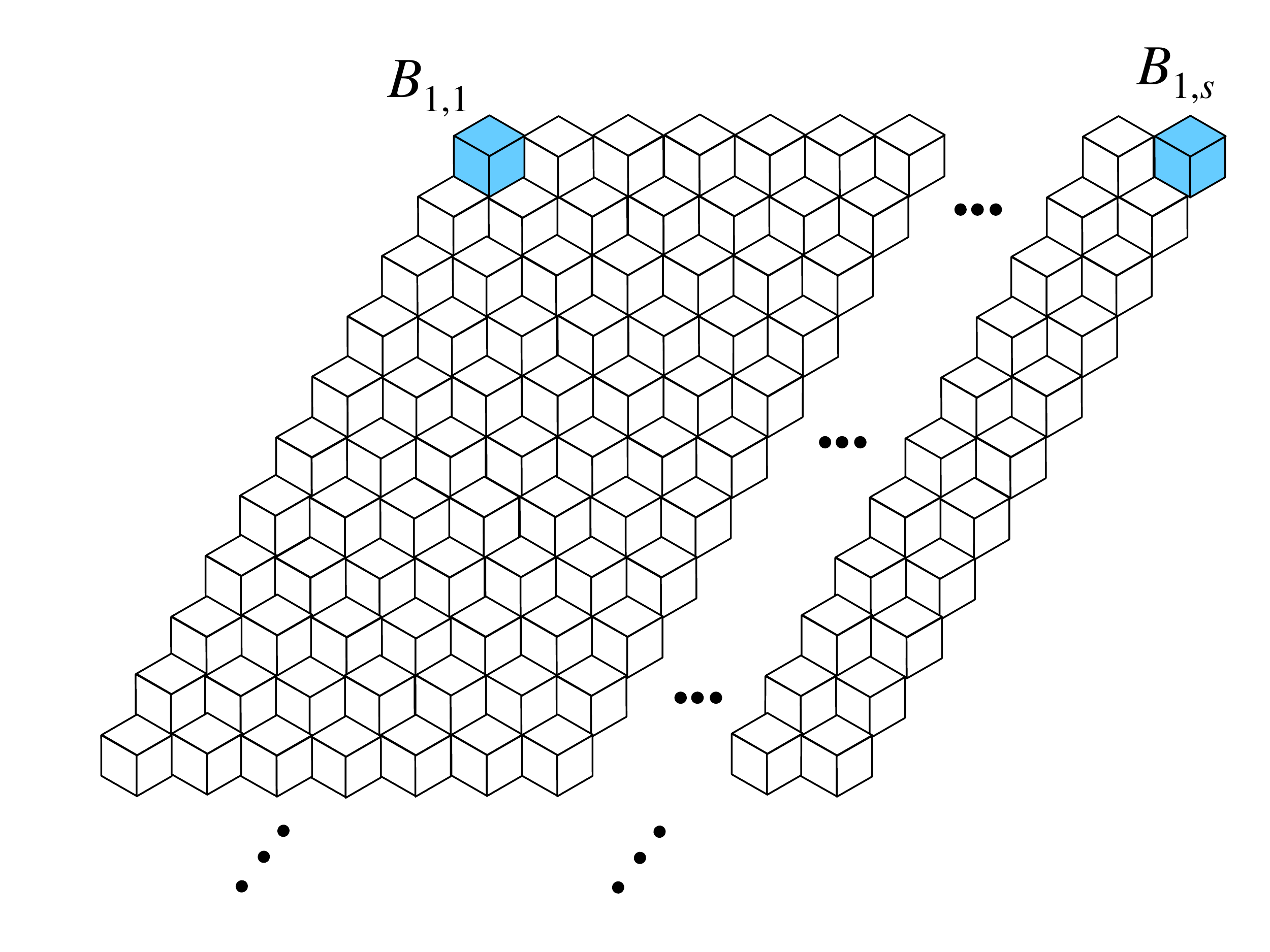} & \includegraphics[width=0.35\textwidth]{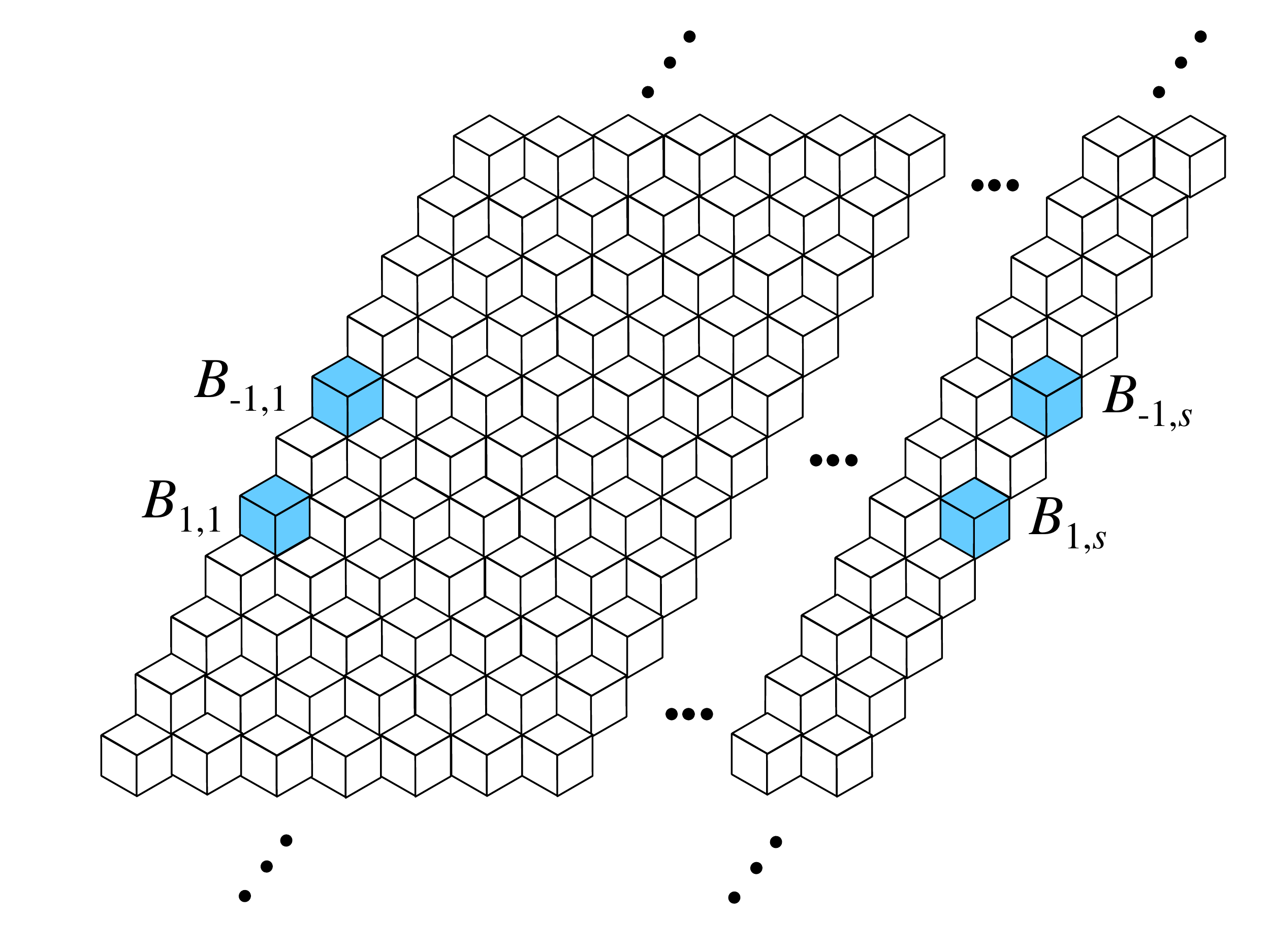} \\
         (a) & (b)
    \end{tabular}
    \centering
    \begin{tabular}{c@{\hskip 5em}c}
     \includegraphics[width=0.35\textwidth]{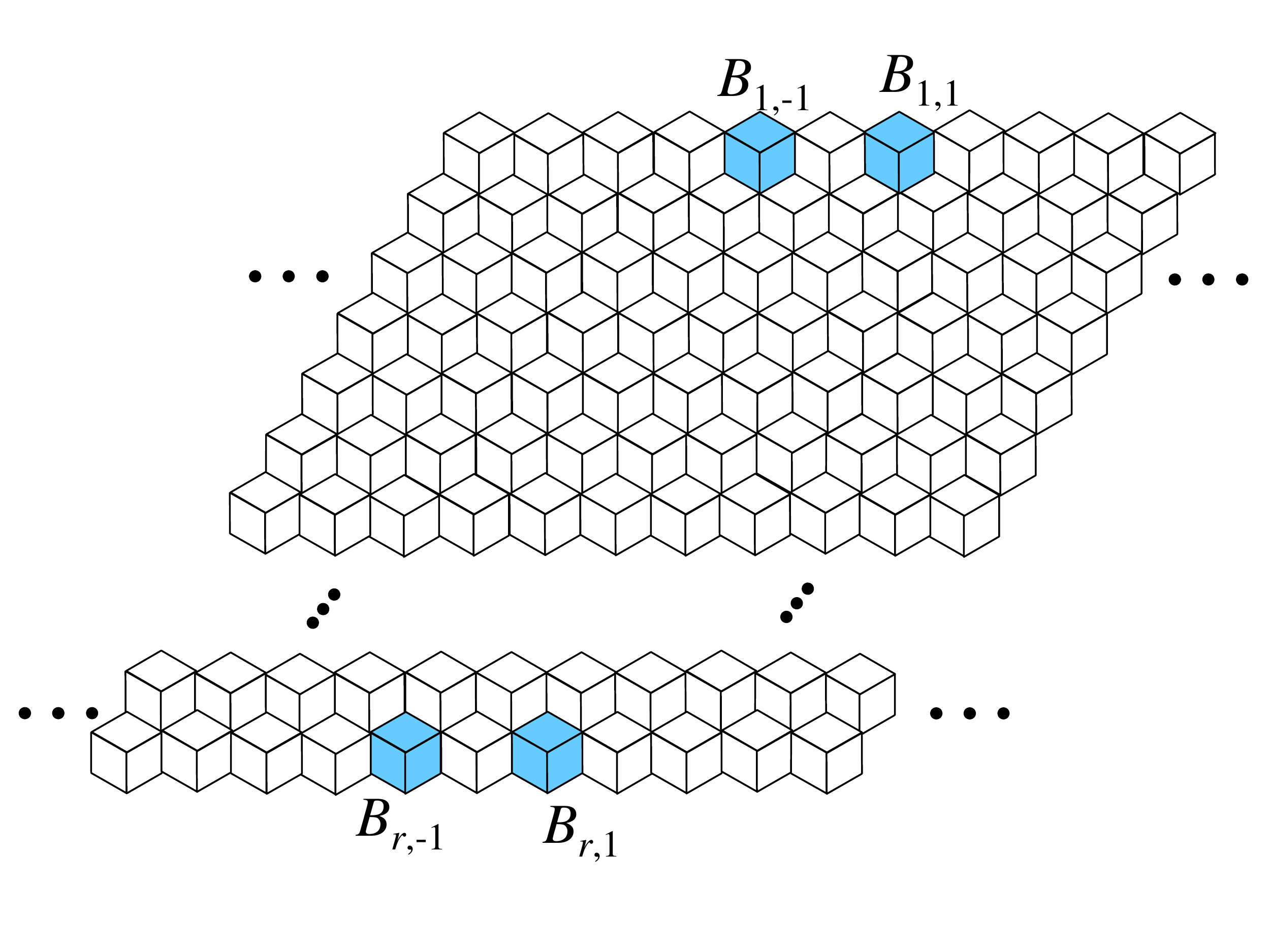} & \includegraphics[width=0.35\textwidth]{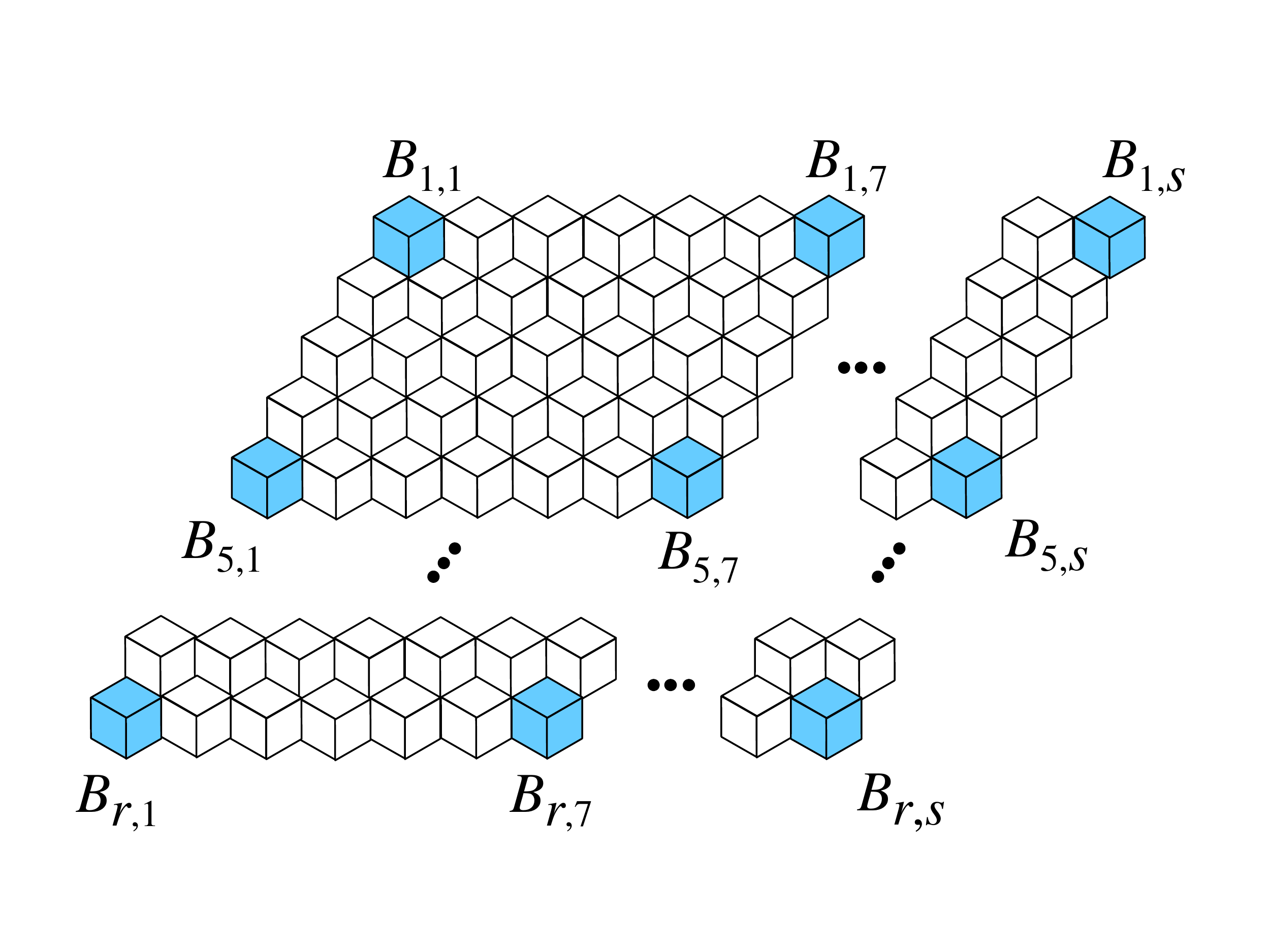} \\
     (c) & (d)
    \end{tabular}
    \caption{Various tumbling block parallelograms}
    \label{fig:tbp}
\end{figure}

Figure \ref{fig:tbp} shows some examples of the tumbling block parallelograms: (a) the one-way infinite $TBP(\infty, s)$,
(b) the two-way infinite $TBP(\infty, s)$, (c) the two-way infinite $TBP(r, \infty)$,
and (d) the tumbling block parallelogram with 5 rows and 7 columns, $TBP(5, 7)$,
within the tumbling block parallelogram with $r$ rows and $s$ columns, $TBP(r, s)$.
Note that $TBP(\infty, s)$ and $TBP(s, \infty)$ are isomorphic.
It can easily be verified that the number of vertices in the
tumbling block parallelogram with $r$ rows and $s$ columns is $|V(TBP(r,s))| = 7 + 5(s-1) + (r-1)[5 + 3(s-1)]	
= 3rs +2r +2s$.  The number of edges is $|E(TBP(r,s))| = 9 + 8(s-1) + (r-1)[8 + 6(s-1)]	 = 6rs +2r +2s - 1$.

\FloatBarrier

\begin{figure}[H]
\centering
% \subfigure[] % caption for subfigure a
% {
%     \label{fig:sub:a-r}
%     \includegraphics[width=170pt]{fig/Figure5a.pdf}
% }
\subfigure[] % caption for subfigure b
{
    \label{fig:sub:b-r}
    \includegraphics[width=180pt]{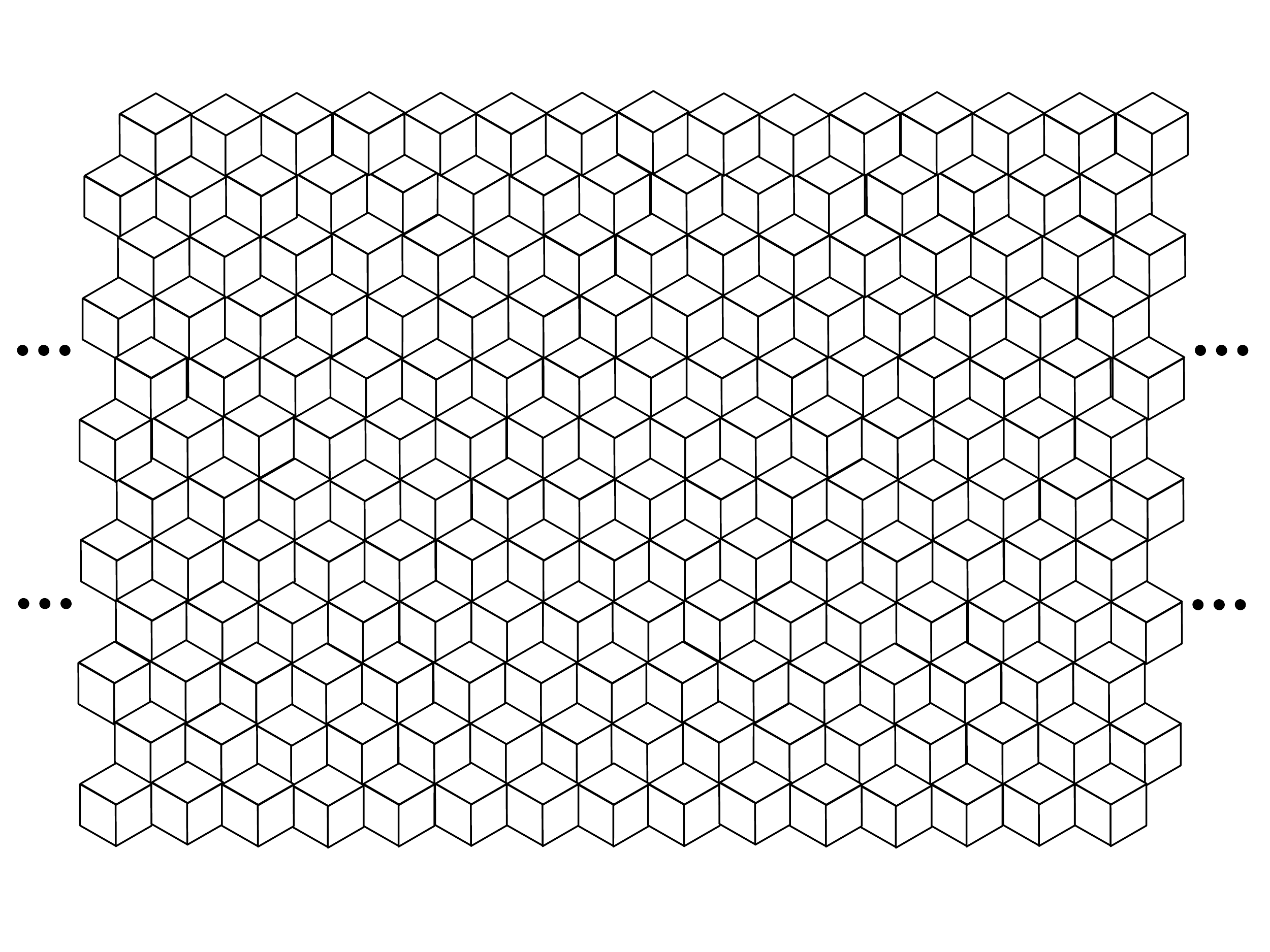}
}
\subfigure[] % caption for subfigure c
{
    \label{fig:sub:c-r}
    \includegraphics[width=180pt]{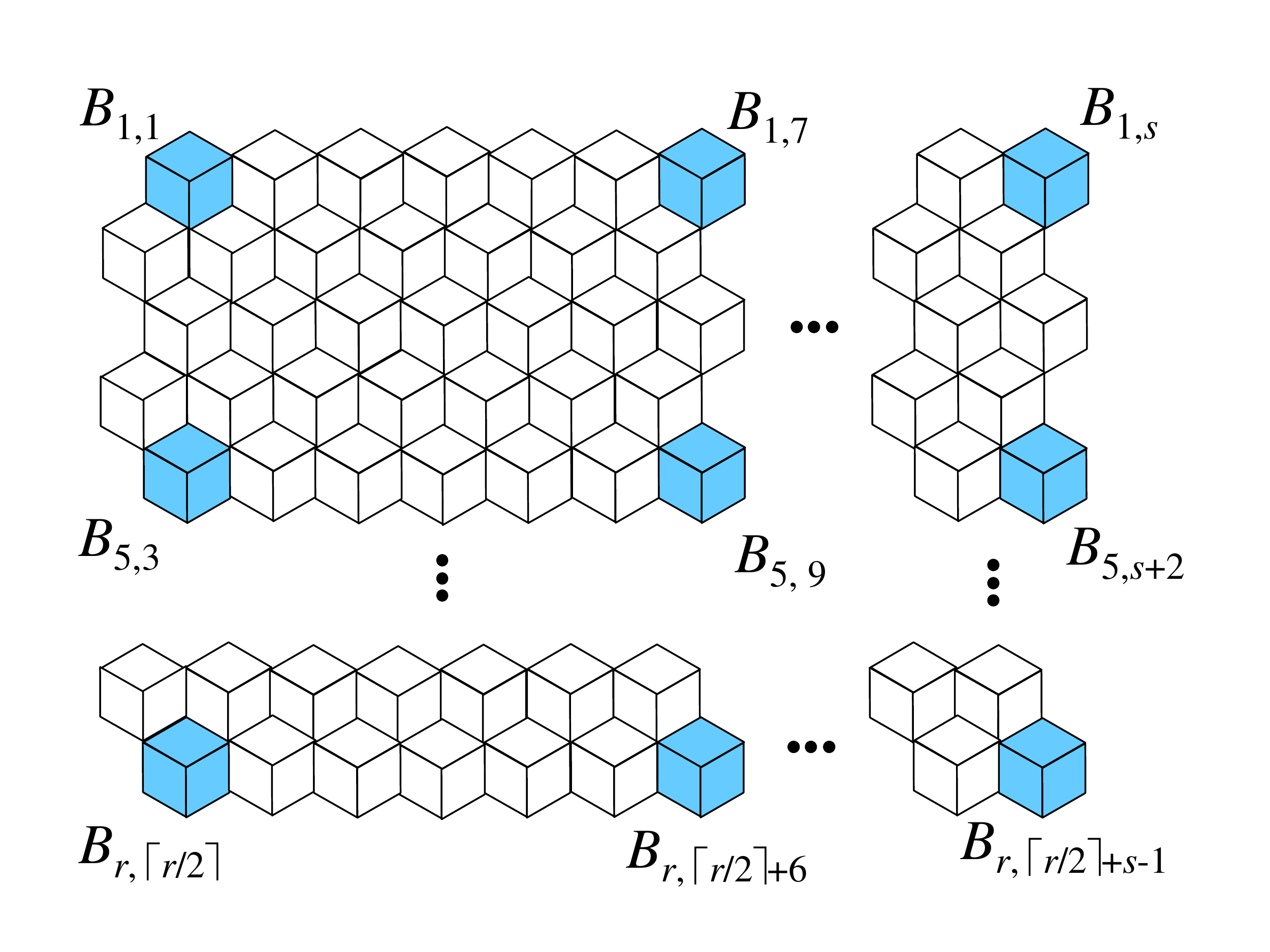}
}
\caption{Various tumbling block rectangles}
\label{fig:tbr} % caption for the whole figure
\end{figure}
\FloatBarrier

Figure \ref{fig:tbr} shows some examples of tumbling block rectangles: (a) the two-way infinite $TBR(\infty, s)$,
(b) the infinite $TBR(r, \infty)$, and (c) the tumbling block rectangle with 5 rows and 7 columns,
$TBR(5, 7)$, within the tumbling block rectangle with $r$ rows and $s$ columns, $TBR(r, s)$.
Note that $|V(TBR(r,s))| = |V(TBP(r,s))|
= 3rs +2r +2s$ and $|E(TBR(r,s))| =|E(TBP(r,s))| = 6rs +2r +2s - 1$.
Observe that $TBR(\infty, s)$ and $TBR(s, \infty)$ are not isomorphic.

Various properties can be considered and various parameters evaluated for these graphs.
We note, for example, that none of these graphs have hamiltonian cycles.  Consider removing
vertex $u_{i,j}$.  Specifically, let $S = {\{u_{i,j}, u_{i,j+1}, u_{i-1,j}, u_{i-1,j-1}, u_{i,j-1}, u_{i+1,j}, u_{i+1,j+1}, u_{i,j+2}, u_{i-1,j+1}, u_{i-2,j}, u_{i-2,j-1},}$

\noindent ${u_{i-2,j-2}, u_{i-1,j-2}, u_{i,j-2}, u_{i+1,j-1}, u_{i+2,j}, u_{i+2,j+1}, u_{i+2,j+2}, u_{i+1,j+2}\}}$.
If we remove these 19 vertices we leave 24 isolated vertices, showing that the graph can not be hamiltonian.
Tumbling block graphs with fewer than four rows or columns can also be seen to be non-hamiltonian.

\section{Graphical parameters for tumbling block graphs}

\subsection{Domination related parameters}
The open neighborhood of vertex $v$, denoted by $N(v)$, is the set of vertices adjacent to $v$, and the degree of $v$ is $deg(v)= |N(v)|$.  The closed neighborhood of $v$ is $N[v] = N(v) \cup \{v\}$.  Vertex $v$ dominates itself and its neighbors, that is, every vertex in $N[v]$, and vertex set $D \subseteq V(G)$ is dominating if every vertex is dominated by at least one $v \in D$, that is, $V(G) = \cup _{v \in D} N[v]$. The domination number $\gamma (G)$ is the minimum cardinality of a dominating set.
\FloatBarrier

\begin{figure}[htp]
\centering
\includegraphics [width = 200pt]{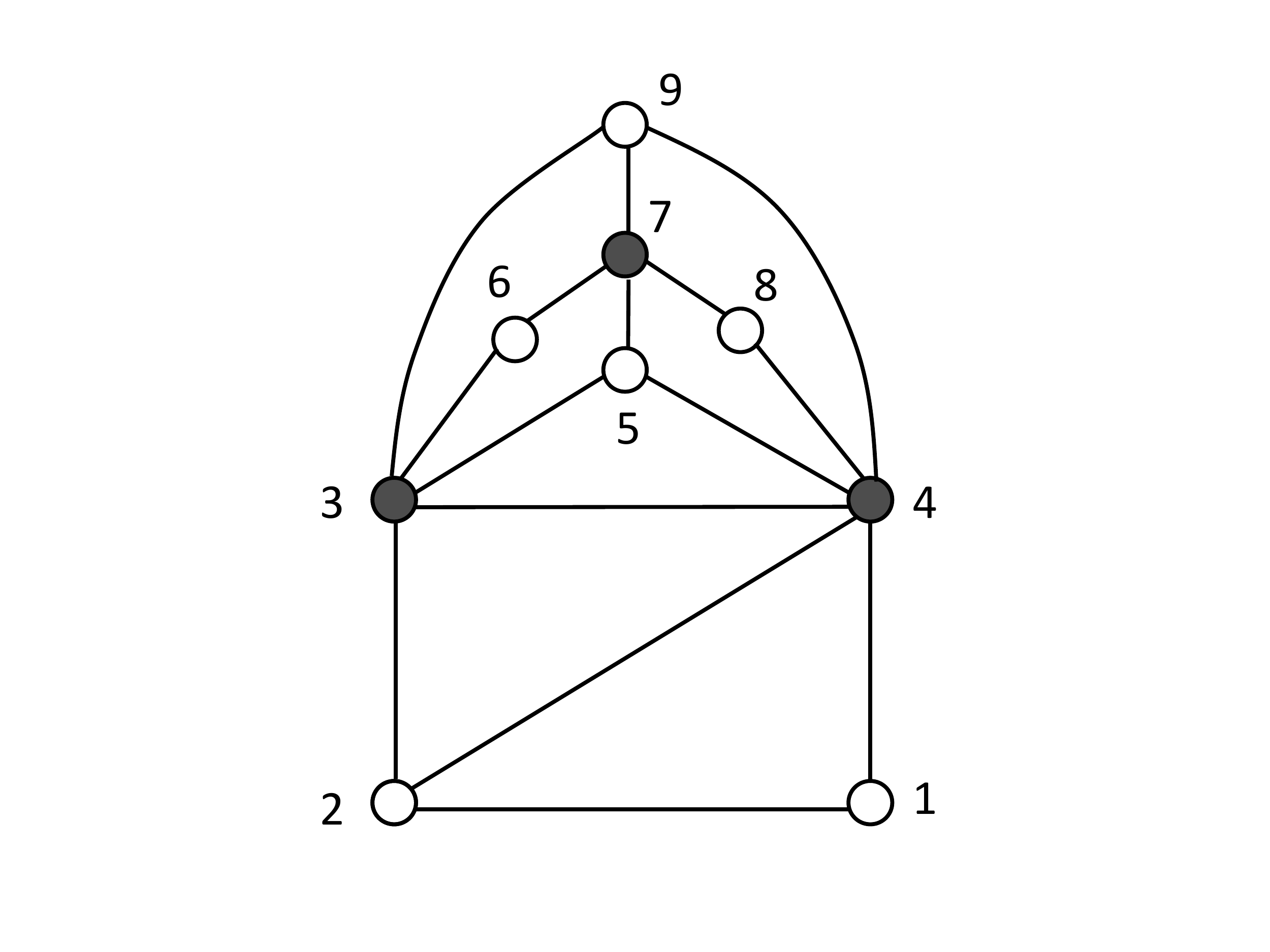}
\caption{Graph $H$} \label{fig:h}
\end{figure}
\FloatBarrier

Percentage parameters for locally-finite, countably infinite graphs
were defined in Slater [19].  For example, for the $\gamma (G)$ parameter we have $\gamma\%(G)$ defined
as the minimum possible percentage of vertices in a dominating set of $G$.
First, we consider $\gamma \%(TB)$ for the infinite tumbling block, $TB$,
and the proof for the next proposition on $\gamma \%(TB)$ is based on \textquotedblleft share" arguments,
as introduced in Slater[19].

If $D$ dominates and $v \in D$,
then the $share$ of $v$ in $D$ is a measure of the amount of domination done by $v$.
If $N[u] \cap D = \{v\}$, then $v$ is a \textit{sole dominator} of $u$ and
$u$ is said to be a \textit{private neighbor} ($PN$) of $v$.
In Figure \ref{fig:h}, vertex 1 is a private neighbor of vertex 4, and vertex 7 is its own $PN$. Because $N[2] \cap \{3,4,7\} = \{3,4\}$, each of vertex 3 and vertex 4 is considered to have a 1/2-share in dominating vertex 2.   If $D$ dominates and $v \in D$, then the $share$ of $v$ in $D$ is defined as $sh(v;D) = \Sigma _{w \in N[v]} 1/ |D \cap N[w]|$.
For example, in graph $H$ of Figure \ref{fig:h} we have $N[3] = \{2,3,4,5,6,9\}$ and
$sh(3; \{3,4,7\})$ = 1/2 + 1/2 + 1/2 + 1/3 + 1/2 + 1/3 = 8/3.
Also, $sh(4; \{3,4,7\})$ = 1 + 1/2 + 1/2 + 1/2 + 1/3 + 1/2 + 1/3 = 11/3,
and $sh(7; \{3,4,7\})$ = 1/3 + 1/2 + 1 + 1/2 + 1/3 = 8/3.
Note that $\Sigma _{v \in D} sh(v; D)$ = $|V(G)|$ = $n$ for any dominating set $D$.

\vspace{5mm}\noindent{\bf Proposition 1}. {For the infinite tumbling block, $TB$, $1/7 < \gamma \%(TB) \leq 1/5 $.}

\begin{figure}[htp]
\centering
\subfigure[] % caption for subfigure a
{
    \label{fig:sub:a-gamma}
    \includegraphics[width=200pt]{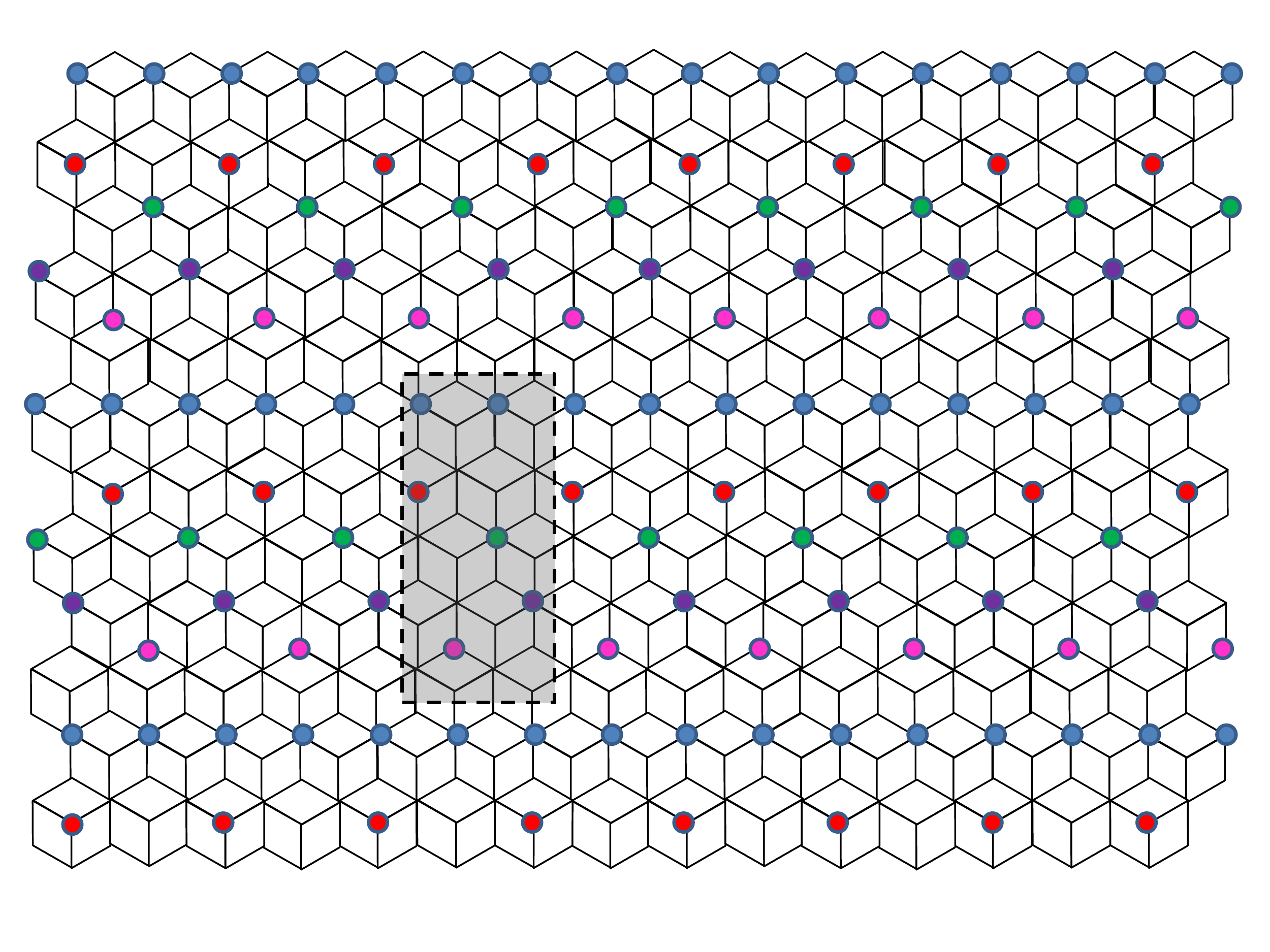}
}
\subfigure[] % caption for subfigure b
{
    \label{fig:sub:b-gamma}
    \includegraphics[width=200pt]{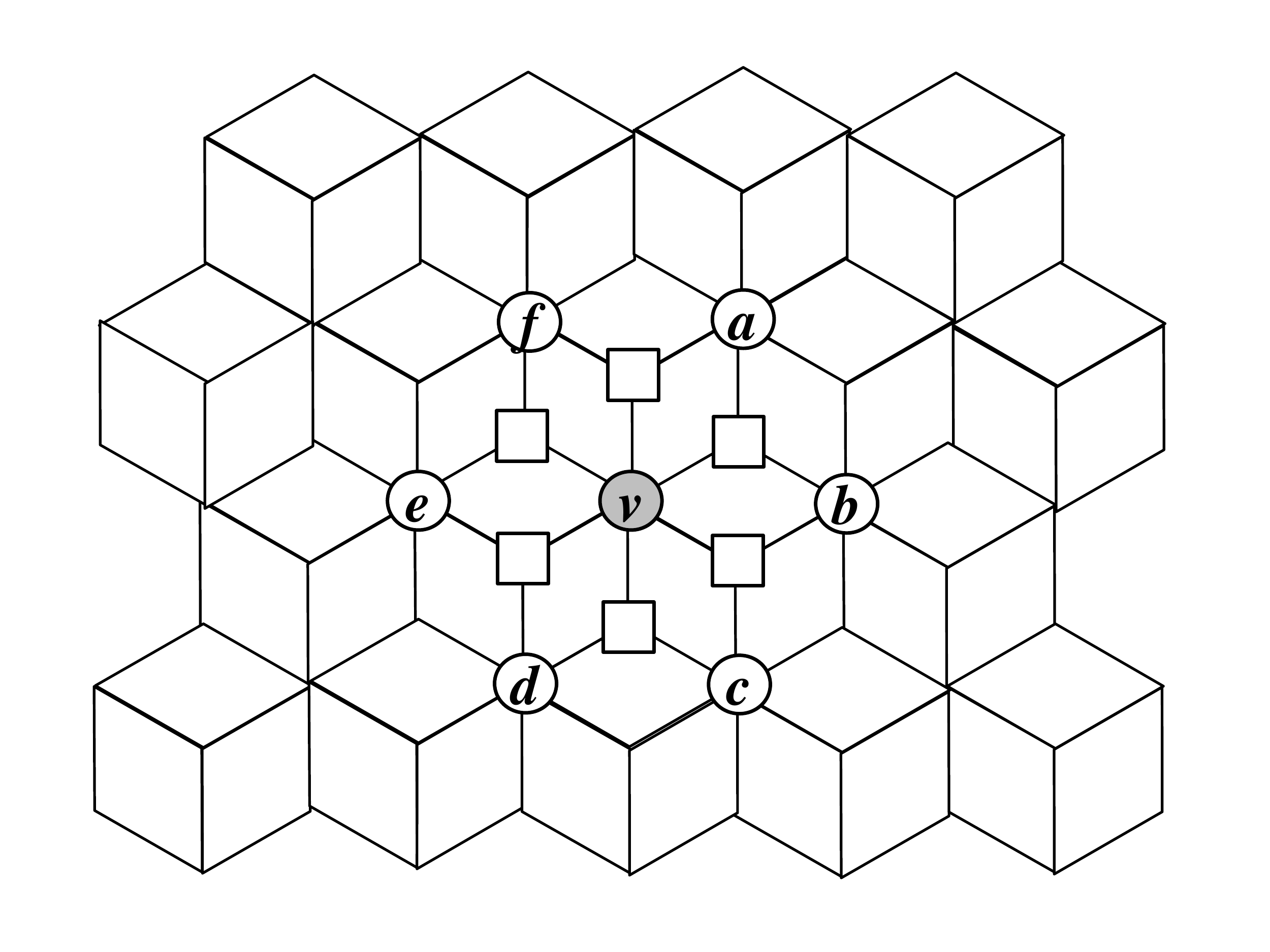}
} \caption{$1/7 < \gamma \%(TB) \leq 1/5 $ }
\label{fig:gamma} % caption for the whole figure
\end{figure}

\begin{proof}
The pattern illustrated in Figure~\ref{fig:sub:a-gamma} shows that we can use 6 vertices out of 30 vertices in a \textit{tiling} of $TB$,
so $\gamma \%(TB) \leq 1/5$.
Note that for a vertex $v$ in a dominating set $D$, $max$  $sh (v;D) \leq 7$ and hence $\gamma \%(TB) \geq 1/7 $.
We could achieve $\gamma \%(TB) = 1/7 $ only when the density of vertices $v$ such that $sh(v;D) = 7$ is one. We will show this is not possible.

Clearly $sh(v;D) = 7$ implies the degree of vertex $v$, $deg (v) = 6$.
Assume $v \in D$ is of degree 6 with $sh (v;D) = 7$.
Let vertex $v$ be such a vertex as shown in Figure ~\ref{fig:sub:b-gamma}.
Then every vertex in $N[v]$ has to be a $PN$ of $v$ shown as square vertices in the figure.
However, the dominating set $D$ must contain at least three vertices of degree 3 in order to dominate the six vertices
 of degree 6 (labeled as $a$, $b$, $c$, $d$, $e$, and $f$ in the figure) that are at distance 2 from vertex $v$,
and the shares of these three vertices are at most 4.
Thus any vertex $v$ with $sh(v;D) =7$ has a vertex at distance two whose share is at most four.
We note that a more complicated share argument can be used to show that
$1/6 \leq \gamma \%(TB) \leq 1/5$, and we believe the upper bound is best possible. \end{proof}

\FloatBarrier

A dominating set $D$ of a graph $G$ is called \textit{efficient} if every vertex is dominated exactly once.
The \textit{efficient domination number} of a graph, denoted by $F(G)$,
is the maximum number of vertices that can be dominated by a set $S$ that dominates each vertex at most once.
For a countably infinite graph $G$ we let $F\%(G)$ denote the maximum possible percentage of vertices
in a set $S$ that dominates each vertex at most once.
A graph $G$ of order $n$  has an efficient dominating set if and only if $F(G) = n$.
A countably infinite graph, $G$ has an efficient dominating set if and only if $F\%(G) = 1$.

\begin{figure}[htp]
\centering
\subfigure[] % caption for subfigure a
{
    \label{fig:sub:a-F}
    \includegraphics[width=200pt]{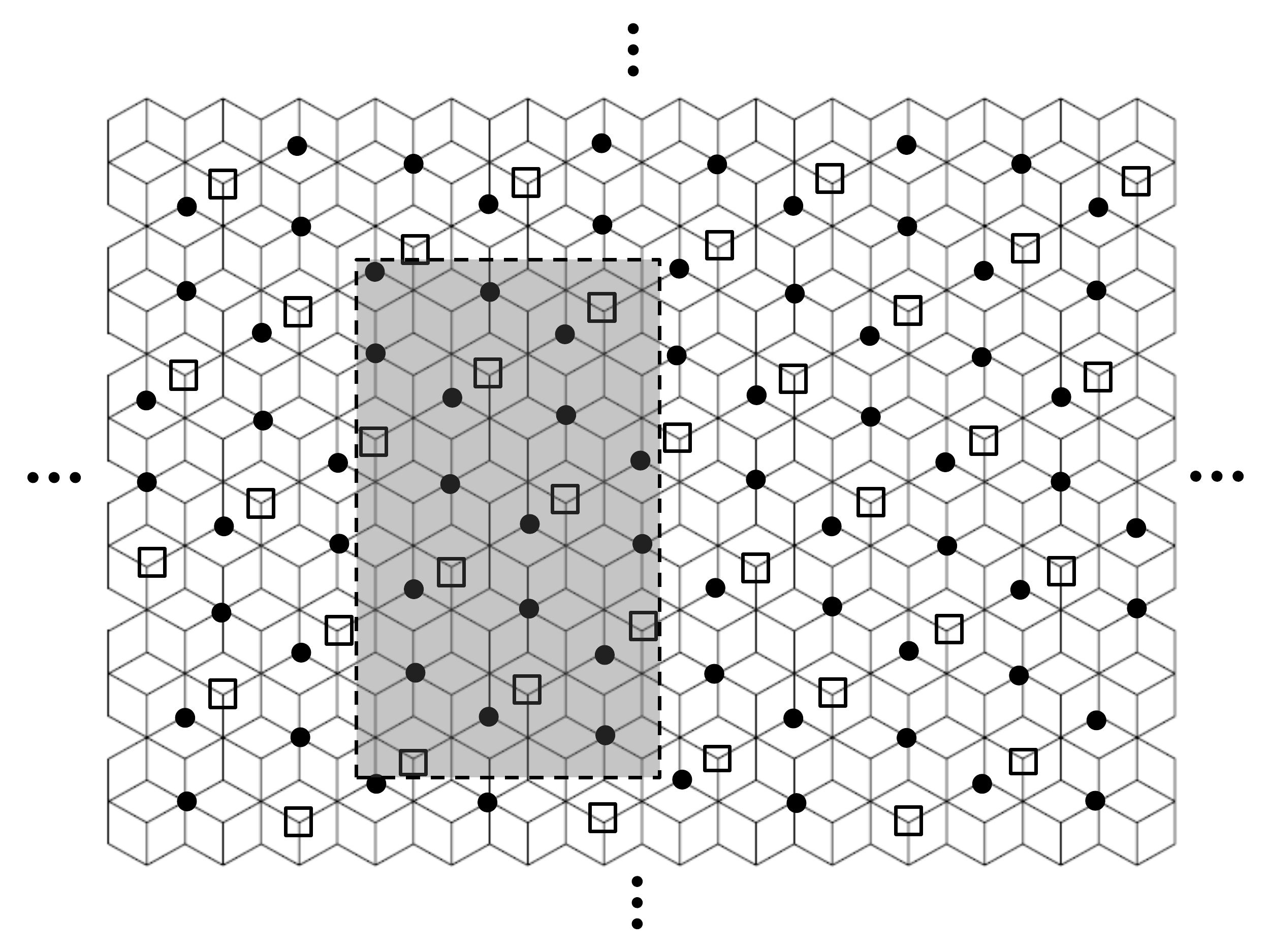}
}
\subfigure[] % caption for subfigure b
{
    \label{fig:sub:b-F}
    \includegraphics[width=200pt]{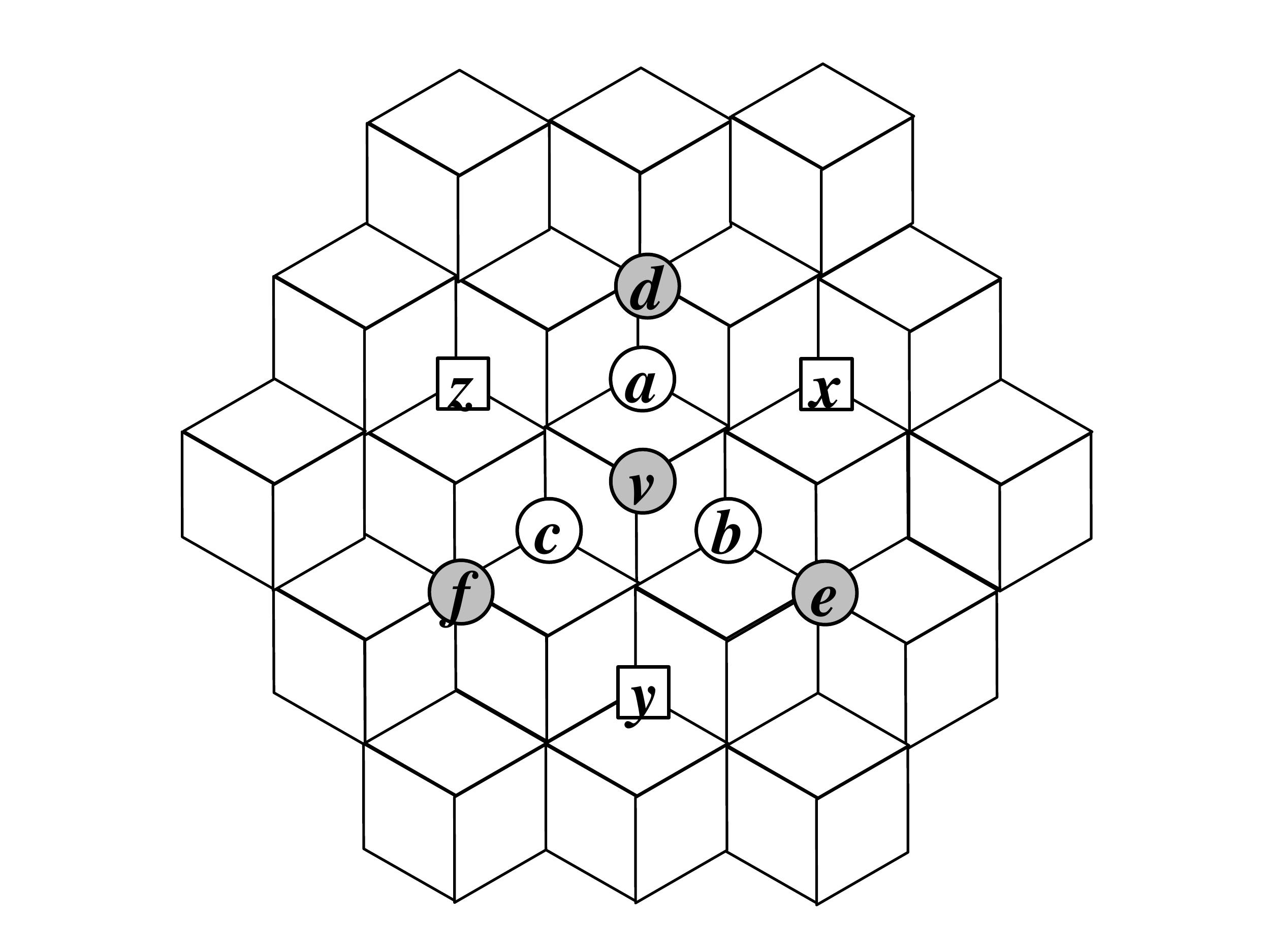}
} \caption{$11/12 \leq F \%(TB) < 1$ }
\label{fig:F} % caption for the whole figure
\end{figure}
\FloatBarrier

\vspace{5mm}\noindent{\bf Proposition 2}. {For the infinite tumbling block, $TB$, $11/12 \leq F \%(TB) < 1$.}

\begin{proof}
The pattern illustrated in Figure~\ref{fig:sub:a-F},
in which $S$ consists of the dark circled vertices and all but the open square vertices are dominated exactly once,
 shows that 88 vertices out of 96 vertices in a \textit{tiling} of $TB$ are efficiently dominated,
so $F \%(TB) \geq 11/12$.  Now we will show $F \%(TB) < 1$.
Assume to the contrary, there is an efficient dominating set, $D$ and
$v \in D$ is a vertex of degree 3 as shown in Figure~\ref{fig:sub:b-F}.
Then, in order to efficiently dominate vertices $a$, $b$, and $c$ we must use vertices $d$, $e$, and $f$.
Now it is impossible to efficiently dominate vertices $x$, $y$, or $z$, so
the density of vertices in $D$ of degree three is zero.
Now it suffices to observe that no two vertices of degree six at distance two can both be in an efficient dominating set $D$.
\end{proof}

Now we investigate $\gamma^{op}\%$ and $F^{op}\%$ parameters for $TB$.
Vertex $v$ \textit{openly dominates} its neighbors, that is, every vertex in $N(v)$ and vertex set $D \subseteq V(G)$
is \textit{open-dominating} (also called \textit{total dominating})
if every vertex is dominated by at least one $v \in D$, that is, $V(G) = \cup _{v \in D} N(v)$.
The open-domination number $\gamma^{op} (G)$ is the minimum cardinality of a open-dominating set.
The efficient open-domination number for a graph $G$, denoted by $F^{op}(G)$, and $F^{op}\%(G)$
for a countably infinite graph are defined similarly to $F(G)$ and $F\%(G)$.

If $D$ open-dominates and $v \in D$, then the \textquotedblleft open-share" of $v$ in $D$ is
a measure of the amount of domination done by $v$. As defined in Seo and Slater [14],
the open-share of a vertex $v$ in open-dominating set $D$, $sh^o(v;D) = \Sigma_{w \in N(v)} 1/ |N(w) \cap  D|$.
For finite graphs $G$ with an open-dominating set $D$ we have $\Sigma_{v \in D} sh^o(v;D) = |V(G)|$
and $|D| \geq  |V(G)| / MAX_{v \in V} sh^o(v;D)$.

\begin{figure}[htp]
\centering
\includegraphics [width = 200pt]{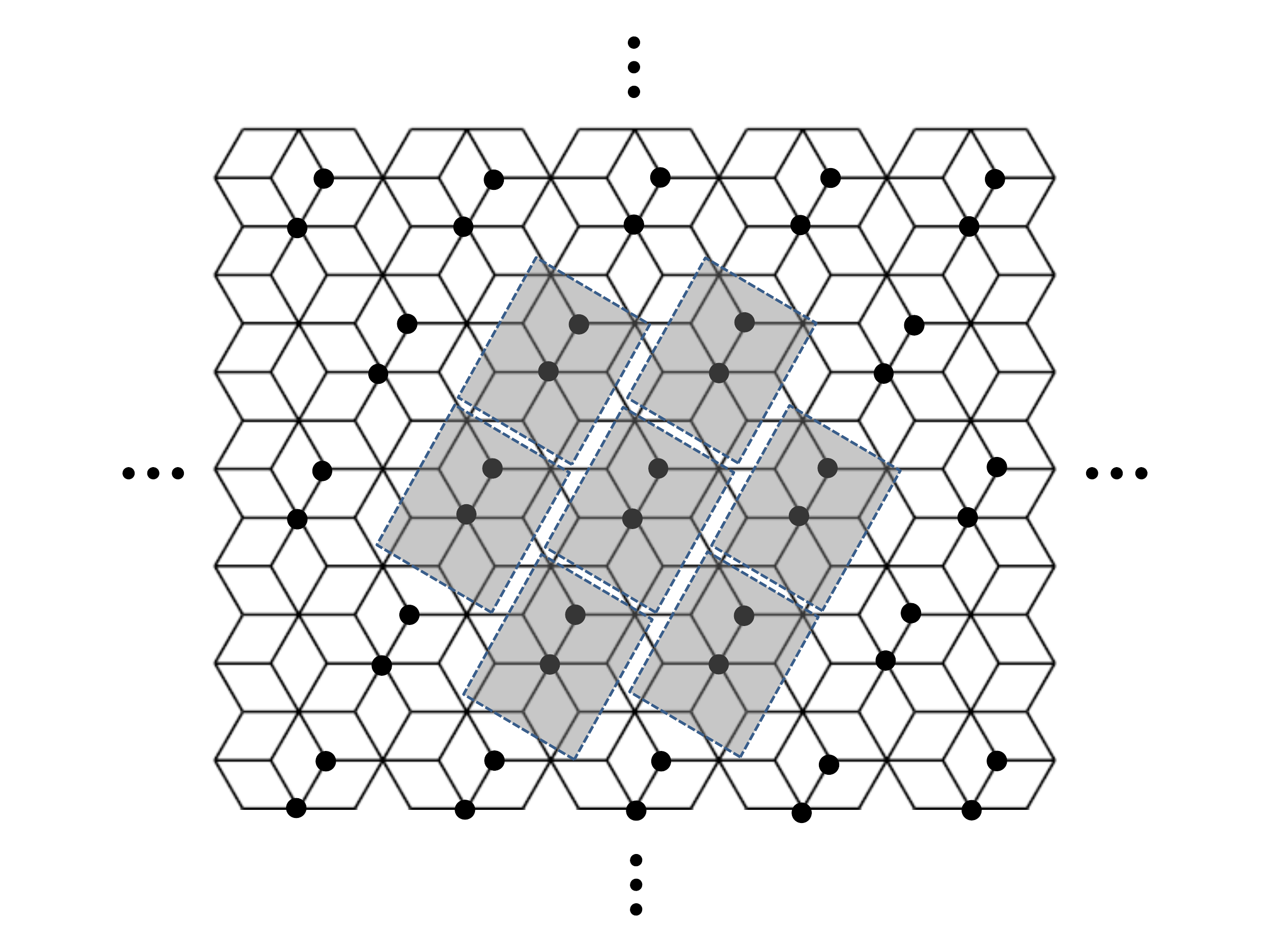}
\caption{$F^{op}\%(TB)=1$ and $\gamma^{op}\%(TB)$ = 2/9}
\label{fig:f-open}
\end{figure}

\vspace{5mm}\noindent{\bf Proposition 3}. {For the infinite tumbling block $TB$, $F^{op}\%(TB)=1$ and $\gamma^{op}\%(TB)$ = 2/9.}

\begin{proof}
The open-dominating set in Figure~\ref{fig:f-open} shows that we can achieve $F^{op}\%(TB)=1$.
This open-dominating set uses 2 vertices out of 9 vertices in a tiling, so $\gamma^{op}\%(TB) \leq 2/9$.
Note that in the infinite tumbling block 2/3 of the vertices are of degree 3, and
the other 1/3 of the vertices are of degree 6.  We observe that $sh^o(v;D) \leq 3$ and $sh^o(v;D) \leq 6$
for a vertex of degree 3 and a vertex of degree 6, respectively.
In order to open-dominate all of the vertices of degree 3 (two-thirds of the entire vertices),
the minimum density of the vertices of $TB$ we need is ((2/3)/6).
Similarly, to open-dominate all of the vertices of degree 6 (one-third of the entire vertices),
 the minimum density of the vertices of $TB$ we need is ((1/3)/3).
Therefore, $\gamma^{op}\%(TB) \geq ((2/3)/6) + ((1/3)/3) = 2/9$, completing the proof.
\end{proof}

\subsection{Distinguishing sets}

For problems involving identifying a malfunctioning processor in a multiprocessor network or
an intruder such as a thief, saboteur or fire in a network modeled facility,
distinguishing sets are of interest.  A collection $S = \{S_1, S_2, ..., S_p\}$ of subsets of $V(G)$
is a distinguishing set for graph $G$ if $\cup_{i=1}^p S_i = V(G)$ and for every pair of distinct vertices $u$ and $v$ in
$V(G)$ some $S_i$ contains exactly one of them.

For locating-dominating sets introduced in Slater [18-20], a detection device at vertex $v$ is assumed
to be able to determine if an intruder is at $v$ or if the intruder is in $N(v)$,
but which vertex location in $N(v)$ can not be determined.
Then $L = \{w_1, w_2, ..., w_j \} \subseteq V(G)$ is a \textit{locating-dominating set} if
$S_1 = \{\{w_1\}, N(w_1), \{w_2\}, N(w_2), ..., \{w_j\}, N(w_j)\}$ is distinguishing.
Identifying codes were introduced in Karpovsky, Chakrabarty and Levitin[11].
For this model a detection device at vertex $v$ can determine if there is an intruder in $N[v]$,
but which vertex location in $N[v]$ can not be determined.
Then $S= \{w_1, w_2, ..., w_j\} \subseteq V(G)$ is an \textit{identifying code}
if $\{N[w_1], N[w_2], ..., N[w_j]\}$ is distinguishing.
When a detection device at $v$ can determine if an intruder is in $N(v)$ but will not report if the intruder is at $v$
we are interested in open-locating-dominating sets as introduced for the $k$-cubes $Q_k$ by Honkala, Laihonen and Ranto [10] and
for all graphs by Seo and Slater [14, 15].  Vertex set $S= \{v_1, w_2, ..., w_j\}$
is an \textit{open-locating-dominating set} if $\{N(w_1), N(w_2), ..., N(w_j)\}$ is distinguishing.
Every graph $G$ has a locating-dominating set; $G$ has an identifying code only when no two vertices have the same closed neighborhood;
and $G$ has an open-locating-dominating set only when no two vertices have the same open neighborhood.
The minimum cardinalities  of a locating-dominating set, an identifying code,
and an open-locating-dominating set are denoted by $LD(G)$, $IC(G)$, and $OLD(G)$, respectively.
Lobstein[12] maintains a bibliography, currently with more than 470 entries, for work on distinguishing sets.

\begin{figure}[htp]
\centering
\subfigure[] % caption for subfigure a
{
    \label{fig:sub:a-ld}
    \includegraphics[width=180pt]{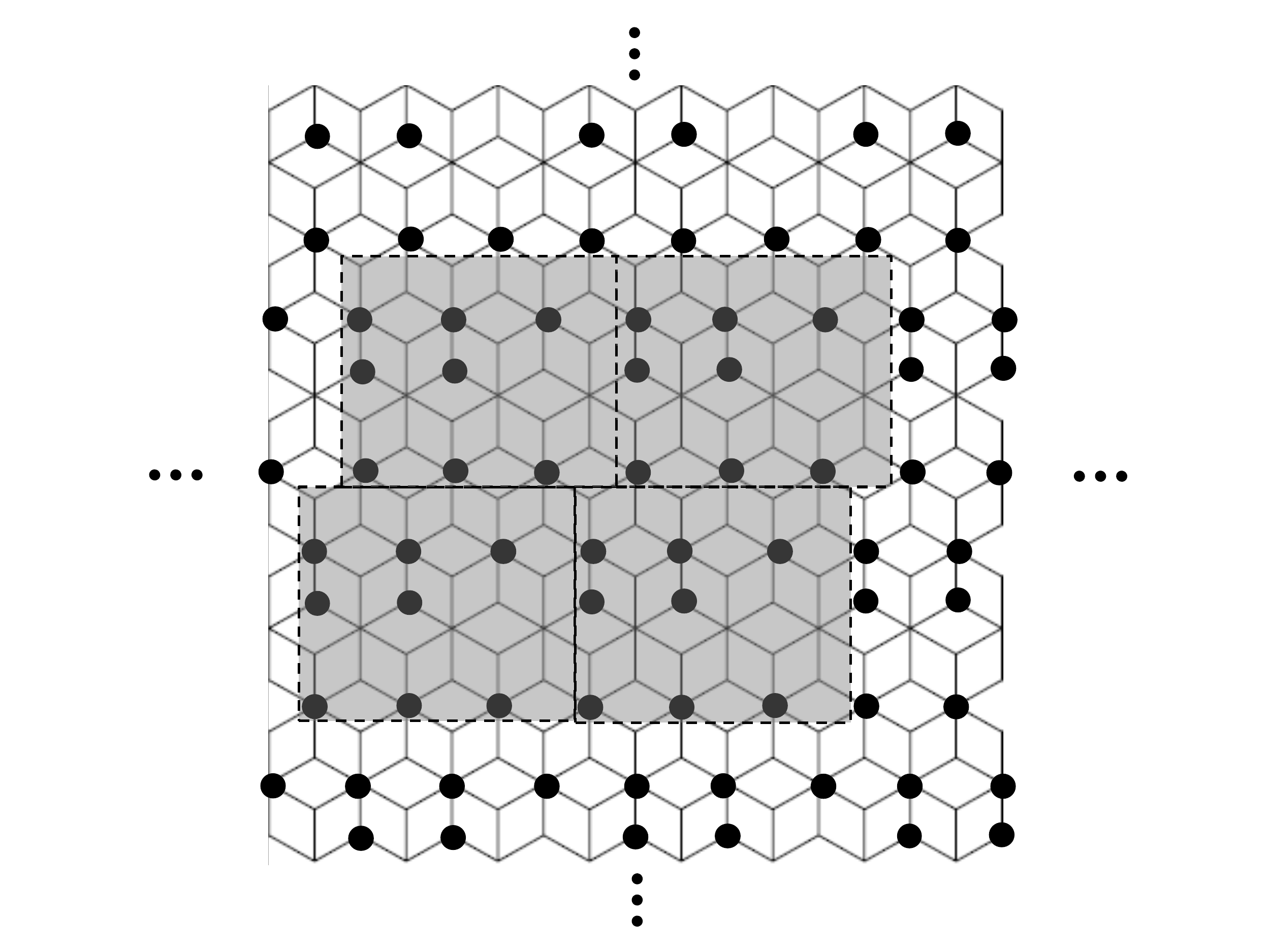}
}
\subfigure[Shares of vertices in $LD$] % caption for subfigure b
{
    \label{fig:sub:b-ld}
    \includegraphics[width=180pt]{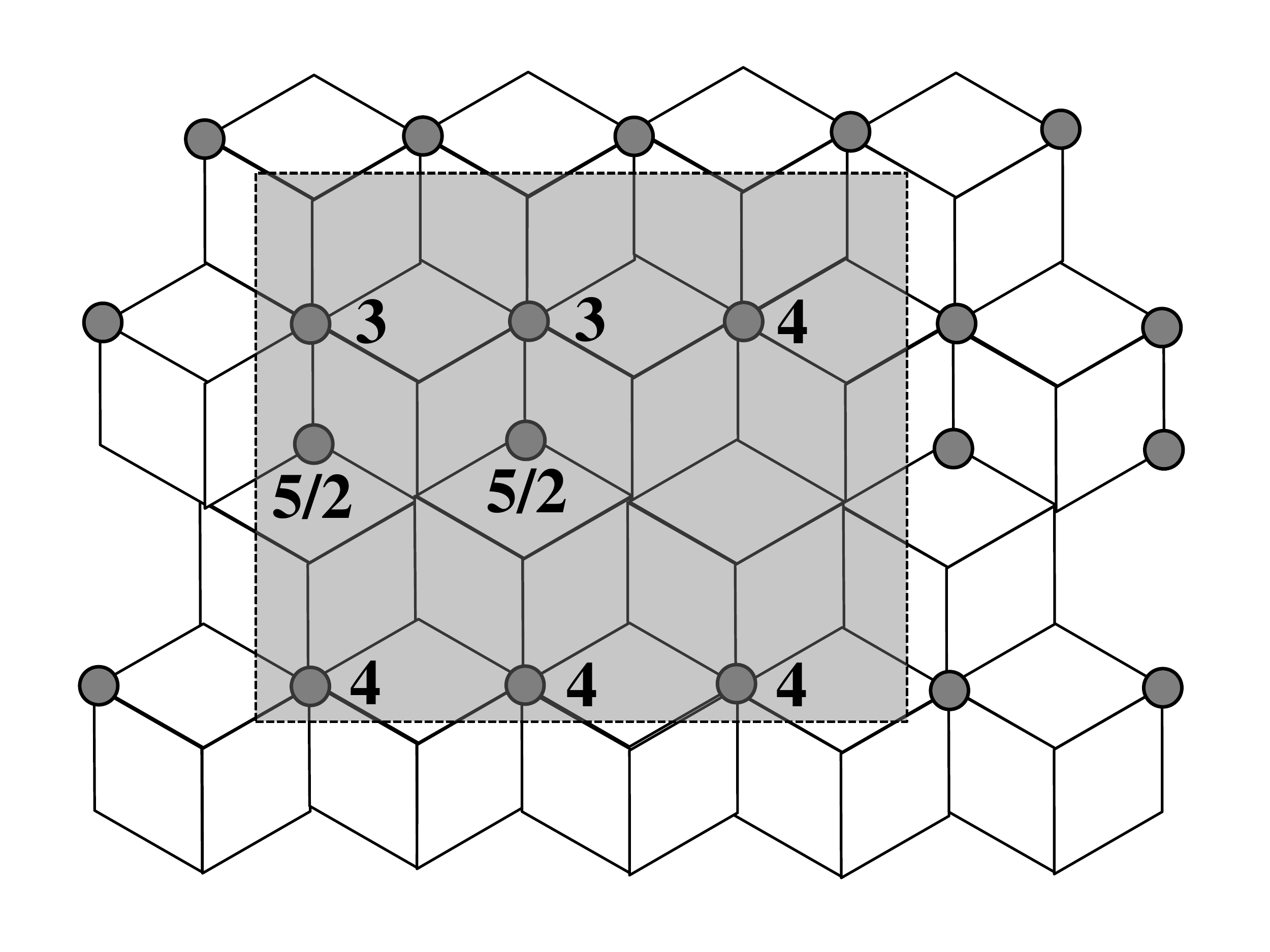}
}
\subfigure[] % caption for subfigure c
{
    \label{fig:sub:c-ld}
    \includegraphics[width=180pt]{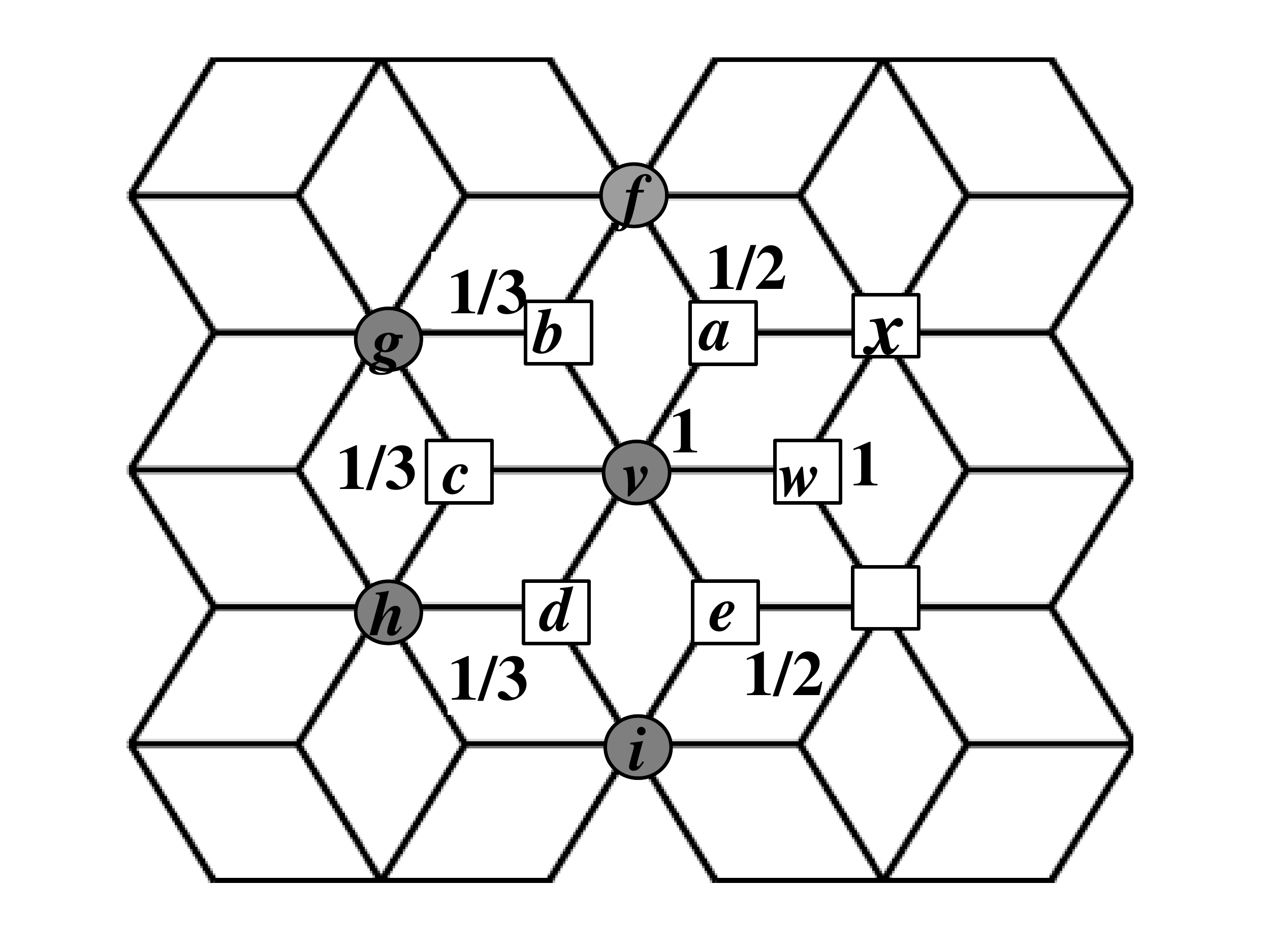}
}
\subfigure[] % caption for subfigure d
{
    \label{fig:sub:d-ld}
    \includegraphics[width=180pt]{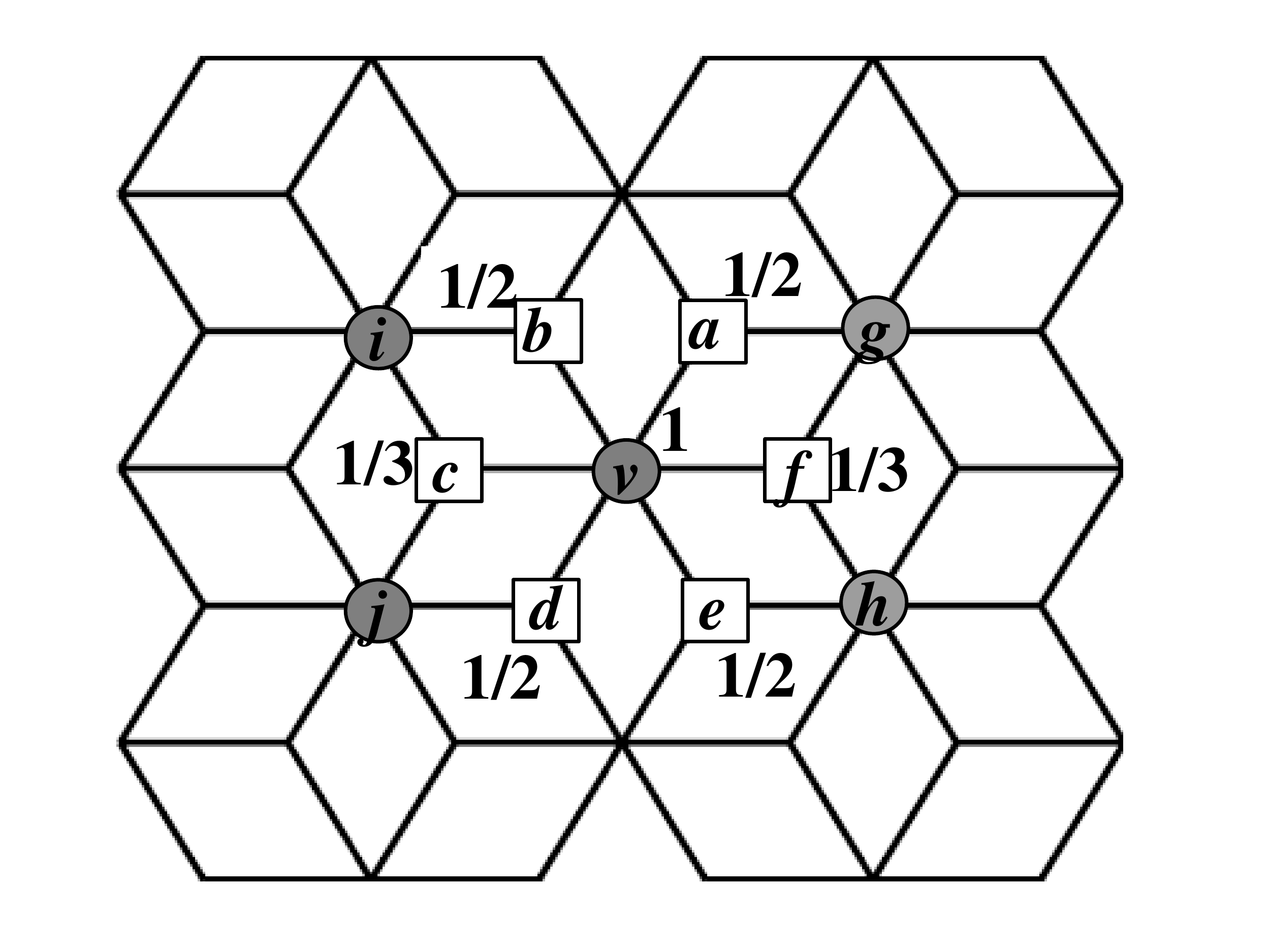}
}
\subfigure[] % caption for subfigure e
{
    \label{fig:sub:e-ld}
    \includegraphics[width=180pt]{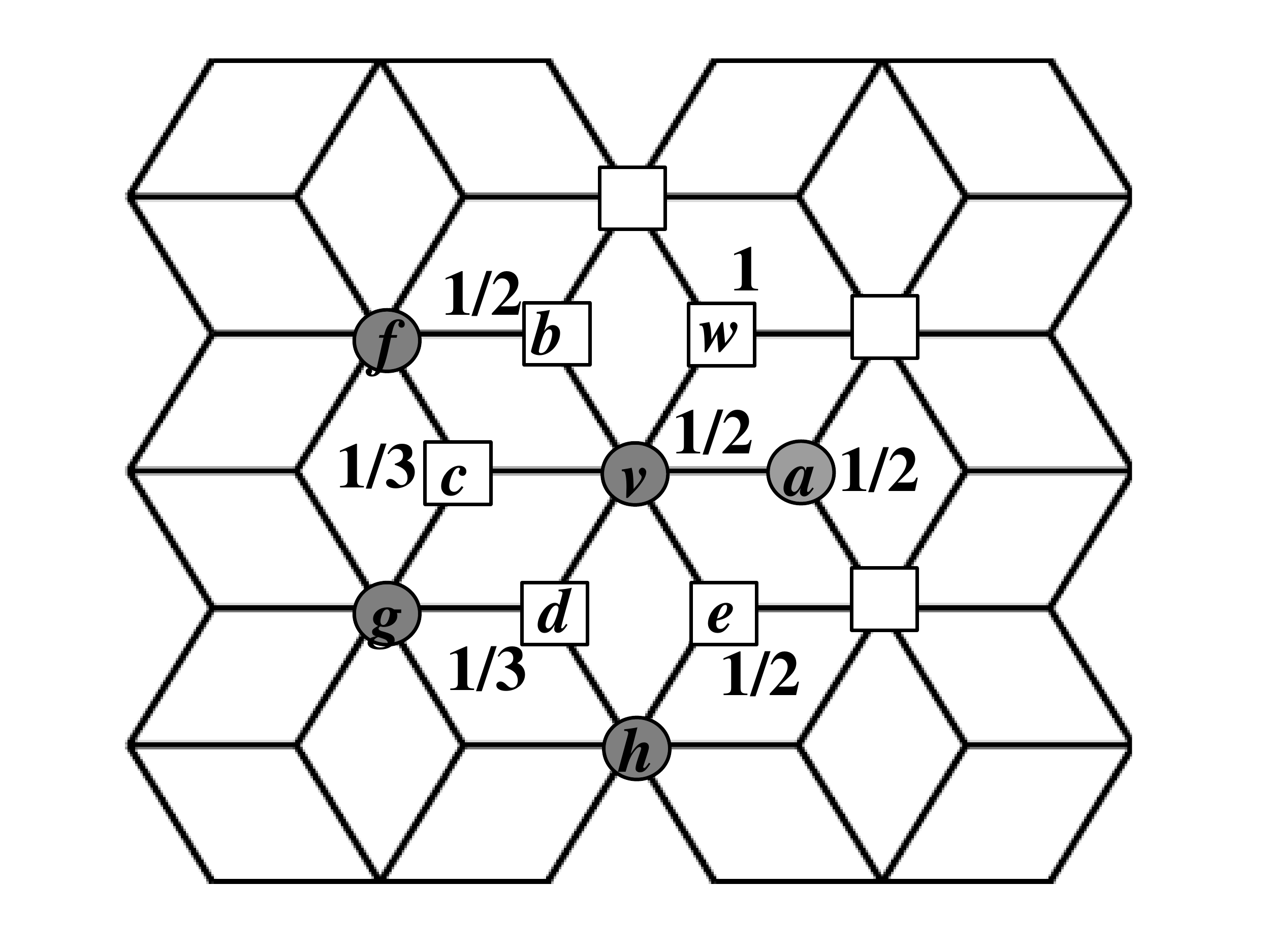}
}
\subfigure[] % caption for subfigure f
{
    \label{fig:sub:f-ld}
    \includegraphics[width=180pt]{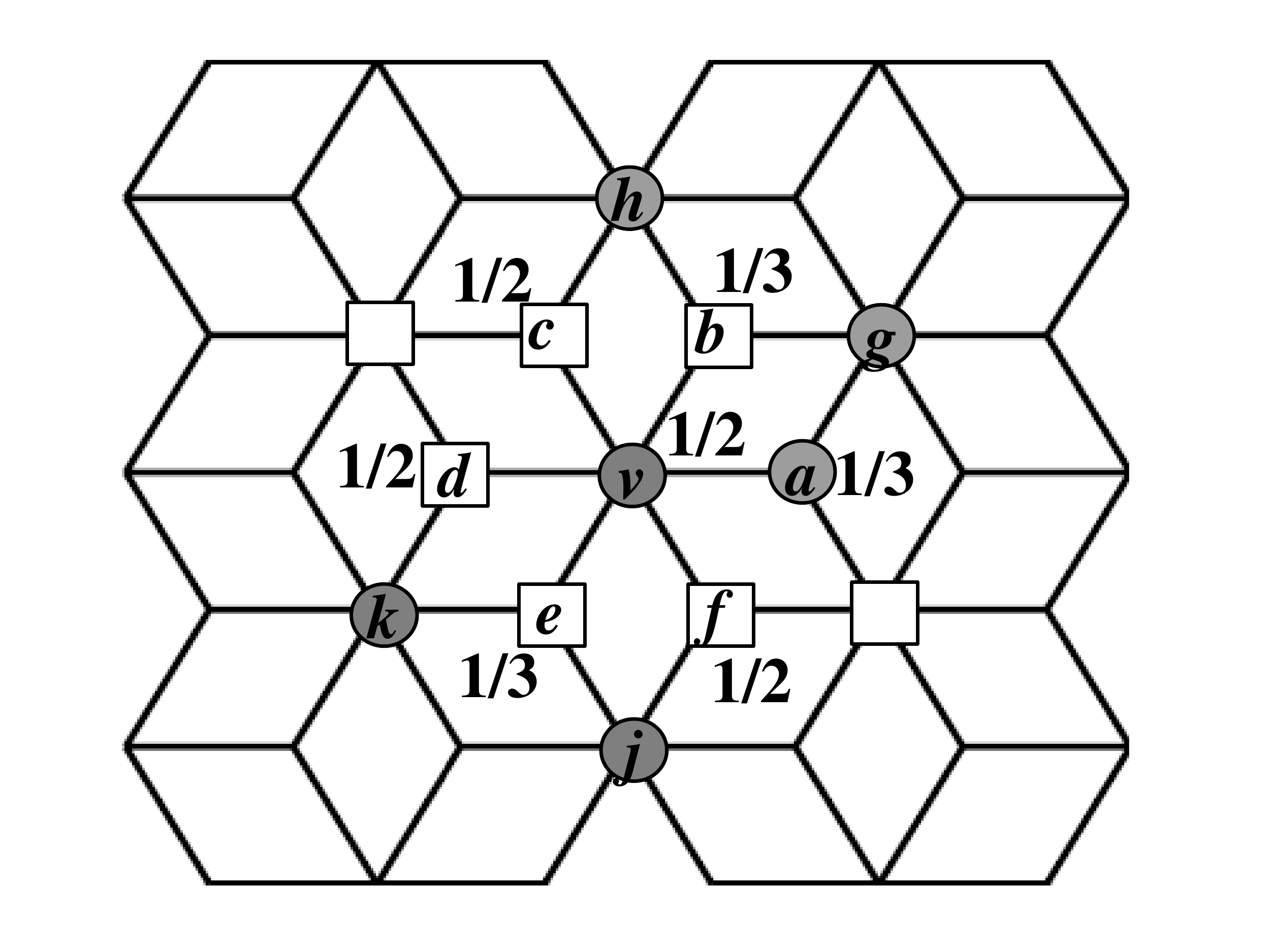}
}
\caption{$1/4 < LD\%(TB) \leq 8/27$ }
\label{fig:ld} % caption for the whole figure
\end{figure}
\FloatBarrier

\vspace{5mm}\noindent{\bf Proposition 4}. {For the infinite tumbling block $TB$, $1/4 < LD\%(TB) \leq 8/27$ .}\label{prop:ld}

\begin{proof}
The pattern illustrated in Figure~\ref{fig:sub:a-ld} shows that eight vertices out of 27 vertices in a \textit{tiling} of $TB$
locate the 27 vertices, so $LD \%(TB) \leq 8/27$.
Another way to look at this is that the total shares of these eight vertices is 27 = 4 + 4 + 4 + 4 + 3 + 3 + 5/2 + 5/2
as shown in Figure~\ref{fig:sub:b-ld}. Now we will show that $LD \%(TB) > 1/4$ using share arguments.
Note that in an $LD$-set each vertex $v$ can have at most one private neighbor in $N(v)$.
It follows that $sh(v) \leq 1 + 1 + 1/2 +1/2 = 3$ and $sh(v) \leq 1 + 1 + 1/2 +1/2 + 1/3 +1/3 + 1/3 = 4$
for $v \in LD$ of degree 3 and degree 6, respectively.  To achieve $LD \%(TB) = 1/4$,
the density of vertices of degree 3 in an $LD$-set must be zero.
Further, to achieve $LD\%(TB) \leq 1/4$ in an $LD$-set $S$ the vertices $v \in S$ with $sh(v;S) = 4$ must have density one.
Such a $v$ must be its own $PN$ and have another $PN$ (as in Figure~\ref{fig:sub:c-ld}).
Note that if $v$ is not its own $PN$ (as in Figure~\ref{fig:sub:d-ld}) or $v$ has no $PN$ (as in Figure~\ref{fig:sub:e-ld}), then $sh(v) \leq 1 + 1/2 + 1/2 +1/2 + 1/3 +1/3 + 1/3 = 11/3$ and if $v$
is not a $PN$ nor has a $PN$ (as in Figure~\ref{fig:sub:f-ld}), $sh(v) \leq 1/2 + 1/2 + 1/2 +1/2 + 1/3 +1/3 + 1/3 = 3$.
Now suppose every vertex $v \in LD$ has share 4 as shown in Figure~\ref{fig:sub:c-ld},
then for each such $v$, there is another vertex $x$ of degree 6 at distance 2 that is not in $N[a] \cap S$, and $x$ must be dominated by a vertex of degree 3 whose maximum share is 3.
Thus the vertices of degree 3 in $S$ would have positive density, which is a contradiction completing the proof.
\end{proof}

\begin{figure}[H]
\centering
\subfigure[] % caption for subfigure a
{
    \label{fig:sub:a-ic}
    \includegraphics[width=180pt]{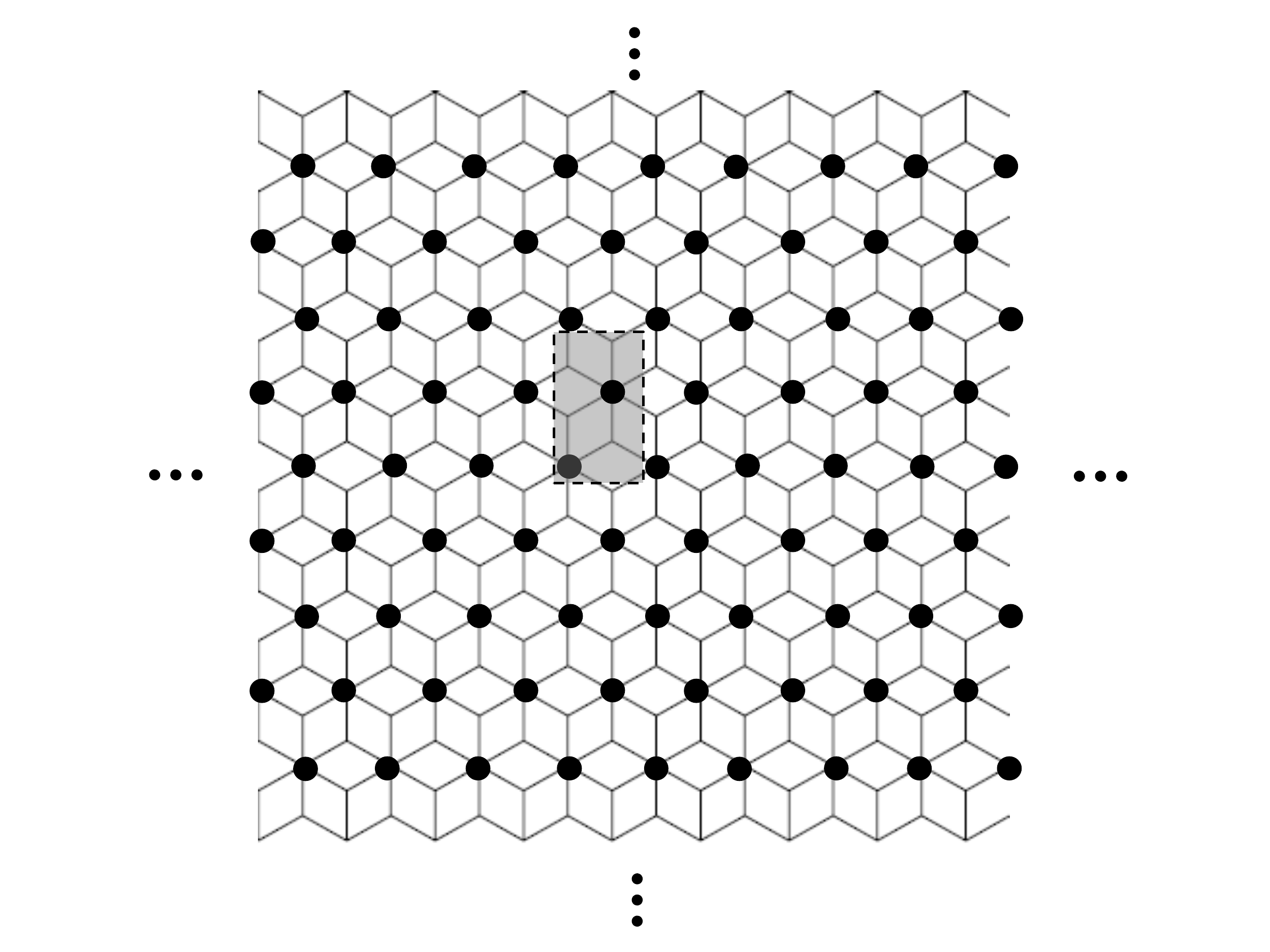}
}
\subfigure[] % caption for subfigure b
{
    \label{fig:sub:b-ic}
    \includegraphics[width=180pt]{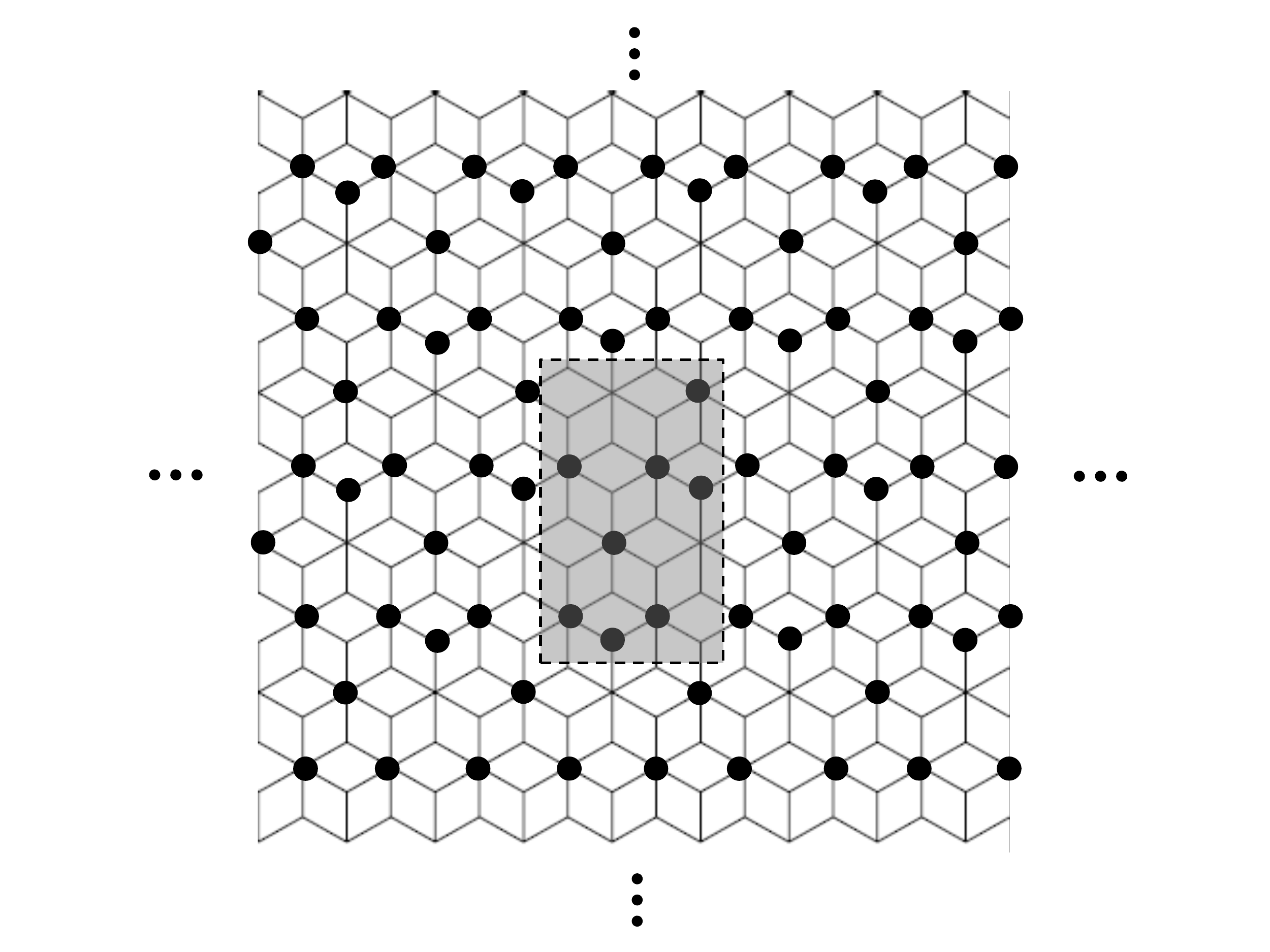}
} \caption{$IC\%(TB) \leq 1/3$ }
\label{fig:ic} % caption for the whole figure
\end{figure}

Using the two different patterns illustrated in Figure~\ref{fig:ic}, we see that $IC\%(TB) \leq 1/3$.  Similar to Proposition 4, we have the following results for $IC\%(TB)$. 

 \vspace{5mm}\noindent{\bf Proposition 5}. {For the infinite tumbling block $TB$, $3/11 \leq IC\%(TB) \leq 1/3$ .}\vspace{5mm}

Next we consider $OLD\%(TB)$ and the proof of Proposition 7 on $OLD\%(TB)$ is based on \textquotedblleft open-share" arguments and the following observation.

\vspace{5mm}\noindent{\bf Observation 6}. {If $D$ is an open-neighborhood locating-dominating set
and $v \in D$, then $sh^o(v;D) \leq   1 + 1/2(deg(v) - 1)$.}

\vspace{5mm}\noindent{\bf Proposition 7}. {For the infinite tumbling block $TB$, $OLD\%(TB) = 7/18$ .}

\begin{figure}[H]
\centering
\subfigure[$OLD\%(TB) \leq 8/27$] % caption for subfigure a
{
    \label{fig:sub:a-old}
    \includegraphics[width=180pt]{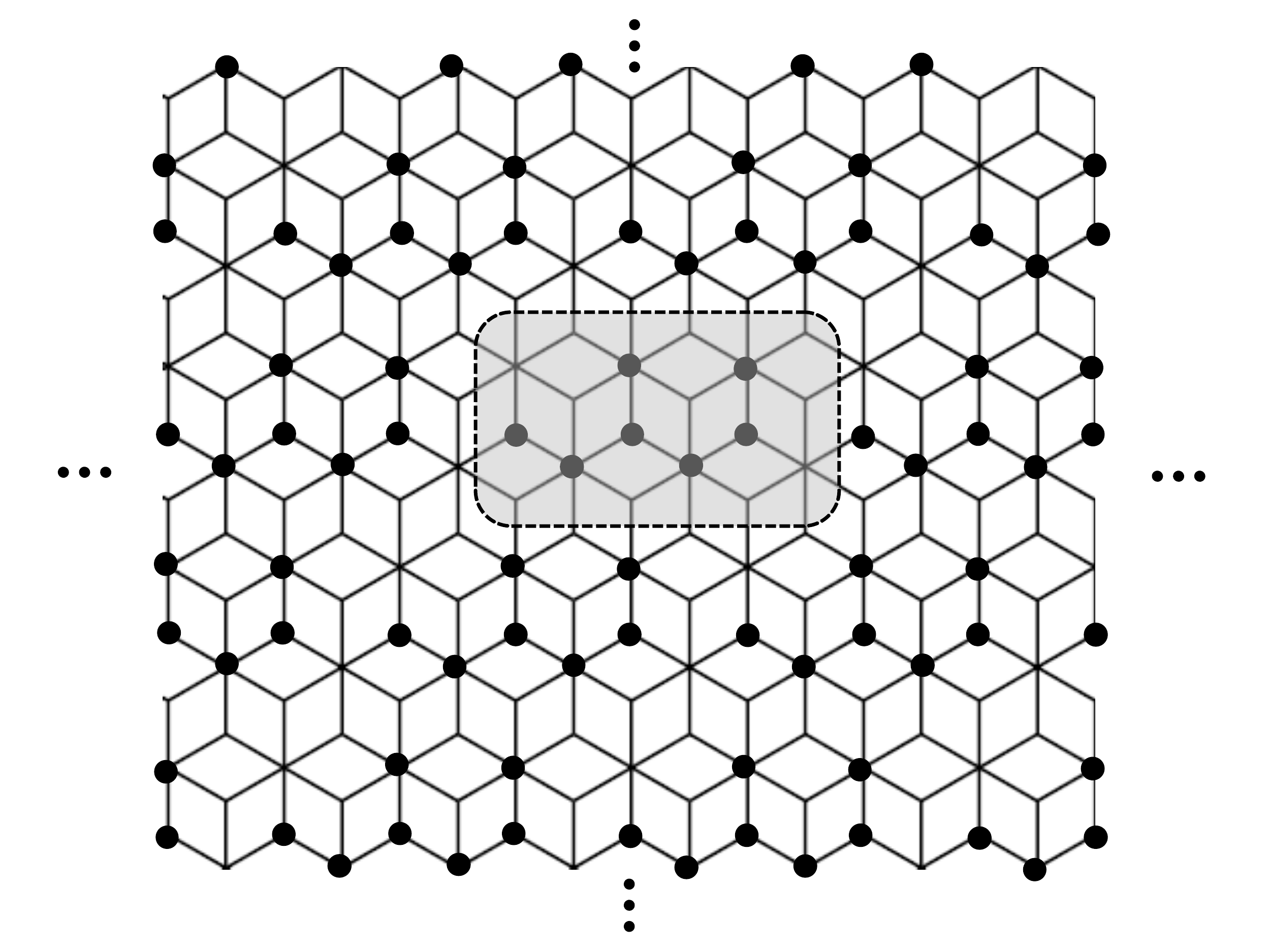}
}
\subfigure[$sh^o (v) \leq 3$] % caption for subfigure b
{
    \label{fig:sub:b-old}
    \includegraphics[width=180pt]{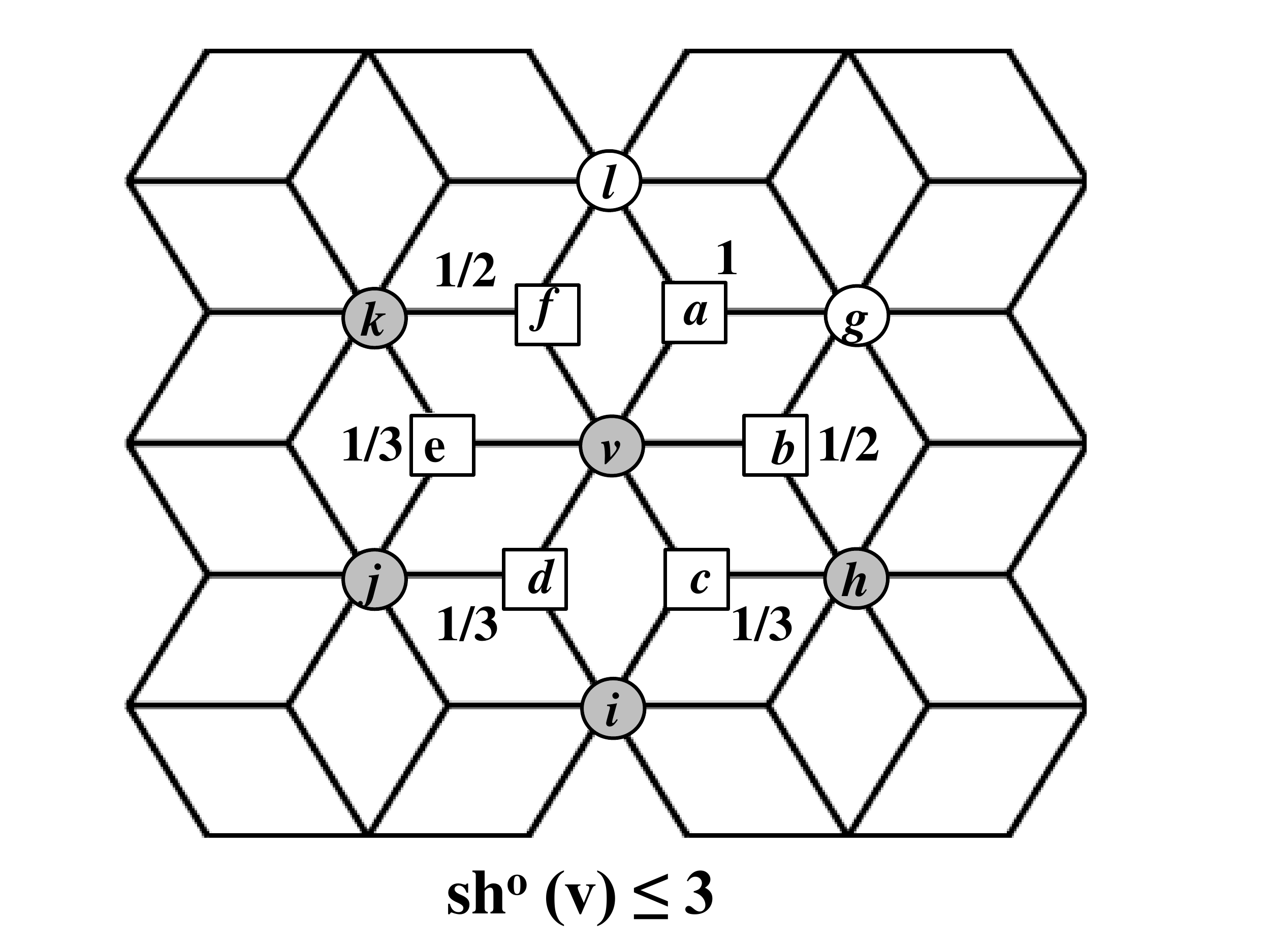}
} \caption{$OLD\%(TB) = 7/18$ }
\label{fig:old} % caption for the whole figure
\end{figure}

\begin{proof}
The pattern illustrated in Figure~\ref{fig:sub:a-old} shows that 7 vertices out of 18 vertices in a \textit{tiling} of $TB$
open-locate the 18 vertices, so $OLD \%(TB) \leq 7/18$. Now we will show that $OLD \%(TB) \geq 7/18$ using
open-share arguments.

Let $S$ be an $OLD(TB)$-set.
By observation 6, for vertex $v$ of degree 3,  we have $sh^{o}(v) \leq 2$.
The one-third of the vertices in $V_2$ that have degree 6 must be open-dominated by the vertices in $V_1$
of degree 3, each with share at most 2.  Hence the density of vertices from $V(TB)$
needed to open-dominate $V_2$ is (1/2)(1/3) = 1/6.

Similarly, $deg(v) = 6$ implies $sh^{o}(v) \leq 7/2$.
However, if $v \in V_2$ has $deg(v) = 6$ and $v$ has a $PN$, say $a$ is a $PN$ of $v$ as in Figure~\ref{fig:sub:b-old}.
Then $N(a) \cap S = \{v\}$, so $g \notin S$ and $N(b) \cap S \neq N(a) \cap S$ implies
that $h \in S$.
Now $N(b) \cap S = \{v, h\} \neq N(c) \cap S$ implies that $i \in S$.
Likewise $\{j, k\} \subseteq S$.
Then $\{v, h, i, j, k\} \subseteq S$ implies that $sh^{o} (v) \leq  1 + 1/2 +1/2 + 1/3 + 1/3 + 1/3 = 3$.
If $v$ has no $PN$, then $sh^{o} (v) \leq  6(1/2) = 3$.
Thus the density of vertices from $V(TB)$ needed to
open-dominate $V_1$ is (1/3)(2/3) = 2/9.
Therefore, $OLD\%(TB) \geq 1/6 + 2/9 = 7/18$, completing the proof.
\end{proof}

\begin{thebibliography}{999}
\bibitem{16}S. Seo and P. Slater, "Open neighborhood locating- domination for grid-like graphs," Bulletin of the Institute of Combinatorics and its Applications, To appear.
\end{thebibliography}
\end{document}